\documentclass{amsart}

\usepackage{latexsym}
\usepackage{amssymb}
\usepackage{amsmath}
\usepackage{fancybox,color}
\usepackage{graphics}
\usepackage{subfigure}

\newtheorem{theorem}{Theorem}[section]

\newtheorem{lemma}[theorem]{Lemma}
\newtheorem{proposition}[theorem]{Proposition}

\newtheorem{definition}[theorem]{Definition}
\newtheorem{remark}[theorem]{Remark}

\def\NN{\mathbb N}
\def\N{\mathbb N}
\def\RR{\mathbb R}
\def\R{\mathbb R}

\def\1{\raisebox{2pt}{\rm{$\chi$}}}
 
\def\a{{\bf a}}
\def\b{{\bf b}}
\def\z{{\bf z}}

\def\zbar{{\bf \overline{z}}}
\def\ubar{\overline{u}}

\def\step{\Delta t}

\def\1{\raisebox{2pt}{\rm{$\chi$}}}

\newcommand{\F}{{\mathcal F}}
\def\v{{\bf v}}

\usepackage{enumerate}

\begin{document}

\title[On a nonlinear flux--limited equation arising in the transport of morphogens]{On a nonlinear flux--limited equation \\arising in the transport of morphogens}

\author{F. Andreu, J. Calvo, J. M. Maz\'{o}n and J. Soler}

\address{Fuensanta Andreu, Jos\'{e} M. Maz\'{o}n
\hfill\break\indent Departament d'An\`{a}lisi Matem\`atica,
Universitat de Val\`encia \hfill\break\indent Valencia, Spain.}
\email{{\tt fuensanta.andreu@uv.es,
mazon@uv.es}}

\address{Juan Calvo, J. Soler  \hfill\break\indent
Departamento de Matem\'{a}tica Aplicada, Universidad de Granada
\hfill\break\indent Granada, Spain.} \email{{\tt juancalvo@ugr.es, jsoler@ugr.es }}



\begin{abstract}
Motivated by a mathematical model for the transport of morphogenes
in biological systems, we study existence and uniqueness of
entropy solutions for a mixed initial-boundary value problem
associated with  a nonlinear flux--limited diffusion system. From
a mathematical point of view the problem behaves more as an
hyperbolic system that a parabolic one.
\end{abstract}

\maketitle


\begin{center}
{\it A Fuensanta, in memoriam. Fuensanta Andreu deceased
26-12-2008}
\end{center}

\section{Introduction}
\setcounter{equation}{0}

The aim of this paper is to analyze the mixed initial-boundary value problem associated with a nonlinear flux--limited reaction--diffusion system
\begin{equation}
\label{DirichletproblemI0} \left\{
\begin{array}{ll}
\displaystyle \frac{\partial u}{\partial t} = \left(\a(u,u_x)
\right)_x - f(t - \tau, u(t,x)) \, u(t,x) + g(t, u(t,x)),
 &
\hspace{0.2cm}\hbox{in \hspace{0.15cm} $]0,T[\times ]0,L[$}\\
\displaystyle
 \\
-  \a(u(t,0),u_x(t,0))  = \beta > 0 \ \ \hbox{and} \ \ u(t,L)=0&
\hspace{0.2cm} \hbox{on \hspace{0.15cm} $t \in  ]0, T[$,}
 \\
\displaystyle
 \\
\displaystyle u(0,x) = u_{0}(x) & \hspace{0.2cm} \hbox{in
\hspace{0.15cm} $x \in ]0,L[$,}
\end{array}
\right.
\end{equation}
being  $$\a(z, \xi):= \nu \frac{ \vert z \vert \xi}{\sqrt{z^2 + \frac{\nu^2}{c^2}\vert
\xi \vert^2}},$$
where the boundary conditions must be interpreted in a weak sense to be precised,
and the functions $f$ and $g$ are nonlinear with respect to $u$ and depend on $u$ through a coupled system of ordinary differential equations. This problem arises in the modelization of the transport of morphogens and the parameter $\tau$ represents a delay in the process of signalling pathway cell internalization.

The nonlinear diffusion equation
\begin{equation}
\label{DirichletproblemI1}
\displaystyle \frac{\partial u}{\partial t} = \left(\a(u,u_x)
\right)_x
\end{equation}
was introduced in different contexts as an alternative
to the linear diffusion equation with the ideas of limiting the flux and reproducing a
system with finite speed of propagation. The flux-limited type equations were motivated
previously in \cite{MM}, but they were firstly deduced by Ph. Rosenau formally, who
proposed three alternative ways to introduce them \cite{Rosenau2}. Also, this equation was
formally derived by  Brenier \cite{Brennier1} by means of Monge-Kantorovich's mass
transport theory  and named {\it relativistic heat equation} after him. As Brenier
pointed out in \cite{Brennier1}, see also \cite{villani}, the relativistic heat equation (\ref{DirichletproblemI1})
is one among the various {\it flux limited diffusion equations} used in the theory of
radiation hydrodynamics \cite{MM}. A general class of flux limited diffusion equations
and the properties of the relativistic heat equation have been studied  in a series of
papers \cite{ACMRelatE,ACMRelat,ACMRelatReg,ACMRelatQ,ACMRelatTM,ACMMRelatDirich}, where
the well--posedness of the Cauchy, the Neumann and the Dirichlet problem for the
relativistic heat equation is proved.

The above discussion on linear diffusion {\it versus} flux--limited diffusion leads to introduce the following change in the classical flux
\begin{equation}
\label{fluxM01} \F = - \nu \nabla u, \qquad \nu > 0,
\end{equation}
associated with the heat equation (or the Fokker-Plank equation)
$$
u_t = \nu \Delta u,
$$
by a flux that saturates as the gradient becomes unbounded. To do that, it was proposed to link $u$
to the flux $\F$ through the velocity $\v$ defined by the relation $\F = u \v.$ Along
with (\ref{fluxM01}) this gives
\begin{equation}\label{FPE1}
\v = - \nu \frac{\nabla u}{u}.
\end{equation}
 According to (\ref{FPE1}), if $\vert \frac{\nabla u}{u} \vert\uparrow \infty$, so will do $\v$.
 However, the
inertia effects impose  a macroscopic upper bound on the allowed free speed, namely, the
acoustic speed or light speed $c$. With this aim, Rosenau proposed to modify (\ref{FPE1})
by taking
\begin{equation}\label{Ros3}
\nu \frac{\nabla u}{u} = \frac{-\v}{\sqrt{1- \frac{\vert \v \vert ^2}{c^2}}}.
\end{equation}
The postulate (\ref{Ros3}) forces $\v$ to stay in the subsonic regime (in the case $c$ is
the acoustic speed). The sonic limit is approached only if $\vert \frac{\nabla u}{u}
\vert \uparrow \infty$. Solving (\ref{Ros3}) for $\v$, we obtain
$$
\F = u \v =  \frac{- u \nabla u}{\sqrt{1+ \left(\frac{\nu \vert
 \nabla u \vert}{c u}\right)^2}}.
$$

As we have mentioned before, the motivation for studying the system (\ref{DirichletproblemI0}) comes from the transport
of morphogenes in biological systems. This is a classical problem since the pioneering work of Turing \cite{T}, Meinhard, Wolpert \cite{wolpert} or  Lander \cite{L2}.  Lander  focussed  the question as a main problem in the understanding of the transport of proteins via signalling pathways: Do morphogen gradients arise by
    diffusion? The relevance of our study is founded on the
analysis of the  Hedgehog (Hh) signaling which has been found to
play multiple roles in development, homeostasis and disease
(reviewed in  \cite{MSR}). In vertebrates the Hh family comprises
three proteins (Sonic, Desert and Indian), which act as secreted,
intercellular factors that affect cell fate, differentiation,
survival, and proliferation in the developing embryo and in many
organs at one time or another. Sonic Hedgehog (Shh) signaling has
also an important role in tumor formation: the deregulation of the
Shh pathway leads to the development of various tumors, including
those in skin, prostate and brain \cite{RA1,RMS,steca}. The idea
is to analyze the morphogenetic patterning of the vertebrate
embryonic neural tube along the dorsoventral (D-V) axis. The
transport of the morphogen Shh along the D-V axis in the neural
tube represents a natural privileged direction for the
description of Shh propagation. Actually,  the system is symmetric
with respect to this axis and this justifies the reduction to one
dimension. The discussion concerning whether the gradient formation of
morphogenes is produced or not by diffusion is a central and
classic topic in developmental biology. This gives a continuous
feedback between mathematical modeling and biological experiments,
see \cite{L2,SS2006,wolpert,GonzalezGaitan}. Recent results in
biology provide some findings that really call into question the
hypothesis of diffussion which has been so often used to model
these phenomena: 1) Concerning the cellular differentiation, the
role of the quantity of morphogen received is as least as relevant
as the time of exposure. With linear diffusion models every point
(cell) of the neural tube receives instantaneously the information
of the morphogen, \cite{DYH,SS2006,PSA}. 2) Morphogenes are transported
in aggregates of several  molecules that also include other
morphogenes or molecules. Then, the typical size of the cluster
aggregates is big (of order 1/10) in comparison with the
extracellular matrix where they are moving { \cite{V}}.
Also their concentration is quite dilute {\cite{guerrero,V,DYH}
}. Therefore, Brownian motion does not seem to
be the more appropriate choice. 3)  In some cases, such as with
the Hh morphogen, it has been proved that in absence of another cell--surface protein, called Ihog, there is
neither propagation nor gradient function of Hh \cite{guerrero}. 4)
There do exist privileged ways/paths of propagation in the
extracellular matrix, a fact that makes the system resemble a
traffic map, more than a linear diffusion system
\cite{vicent,guerrero}.

In this setting, the present paper tries to give some insight on this biological
 problem where
the model here studied is a first step towards a complete model consisting in
$$
   \frac{\partial u(t,x)}{\partial t}= \a(u(t,x),u(t,x)_x)_x - f(t -\tau, u(t,x)) \, u(t,x) + g(t, u(t,x)),
$$
where $f$ stands for the concentration of transmembrane receptor in the cells,
 $g$ represents the concentration of the complex  binding the morphogen to
 the receptor, and where the dependence on $u$ is given through a coupling with a system of seven ODE's
 modeling the rates of change of the concentrations of the proteins participating in the signaling
  pathway coming from the  biochemical  cascade  inside
the cells, see \cite{ArSSV}. In  that work it was also proved that numerical evidence fully agrees with the experiments from a quantitative as well as qualitative (propagation of fronts instead of linear diffusion behaviour) point of view, see \cite{guerrero,shu}.

In addition to the biological or physical motivations, the mathematical analysis of this equation
poses several difficulties, making even more interesting its study, such as the existence and evolution
 of fronts as well as  the  study of its  finite speed of propagation, the related lack of regularity
 and the set--up of an appropriate functional framework to give a meaning to the differential operator
 and the boundary conditions.  In fact, this flux--limited equation provides a behaviour more
  related to hyperbolic systems than to usual diffusive (Fokker--Plack) systems.
  To deal with these mathematical problems we need to combine and extend the applicability of  different
   techniques stemming from parabolic and hyperbolic contexts such as Crandall-Liggett's theorem, Minty-Browder's technique, the concept of entropy solution, and
the  method of doubling variables due to S. Kruzhkov.

This paper deals with a preliminary study of (\ref{DirichletproblemI0}) consisting in  the analysis of the following system
\begin{equation}
\label{DirichletproblemI2} \left\{
\begin{array}{ll}
\displaystyle \frac{\partial u}{\partial t} = \left(\a(u,u_x)
\right)_x ,
 &
\hspace{0.2cm}\hbox{in \hspace{0.15cm} $]0,T[\times ]0,L[$,}\\
\displaystyle
 \\
-  \a(u(t,0),u_x(t,0))  = \beta > 0 \ \ \hbox{and} \ \ u(t,L)=0,&
\hspace{0.2cm} \hbox{on \hspace{0.15cm} $t \in  ]0, T[$,}
 \\
\displaystyle
 \\
\displaystyle u(0,x) = u_{0}(x), & \hspace{0.2cm} \hbox{in
\hspace{0.15cm} $x \in ]0,L[$,}
\end{array}
\right.
\end{equation}
where the boundary conditions must be considered in a weak sense. Our main result is
\begin{theorem} \label{mainresultado} For any initial datum
  $0 \leq u_0 \in   L^{\infty}(]0,L[)$,  there exists a unique
  bounded entropy solution $u$ of
 (\ref{DirichletproblemI2}) in $Q_T = ]0,T[\times ]0,L[$ for every $T
> 0$ such that $u(0) = u_0$.
\end{theorem}

The paper is structured as follows. In the next section we introduce all the tools needed
to develop the theory: a suitable integration by parts formula, lower semi--continuity
results and a functional calculus, in order to be able to give a sense to the
differential operator. In Section 3 we discuss the associated elliptic problem: we define
what a solution is, and then we prove existence and uniqueness of such a solution. Next,
this material is used to define an accretive operator and construct a nonlinear
semigroup, which accounts for solving (\ref{DirichletproblemI2}) in a mild
sense; all this is the content of Section 4. In Section 5 we prove that the mild solution
previously constructed can be  characterized in more operative terms, as a so--called
entropy solution --a concept which is also introduced in this section--, and we prove a
comparison criterium which in particular  entails uniqueness of entropy solutions, thus
proving Theorem \ref{mainresultado}.

\section{Preliminaries}

\subsection{BV functions and integration by parts}

For bounded variation function of one variable we follow \cite{Ambrosio}. Let $I \subset
\R$ an interval, we say that a function $u \in L^1(I)$ is of bounded variation if its
distributional derivative $Du$ is a Radon measure on $I$ with bounded total variation
$\vert Du \vert (I) < + \infty$.
 We denote
by $BV(I)$ the space of all function of bounded variation in $I$. It is well know (see
\cite{Ambrosio}) that given $u \in BV(I)$ there exists $\overline{u}$ in the equivalence
class of $u$, called a good representative of $u$ with the following properties. If $J_u$
is the set of atoms of $Du$, i.e., $x \in J_u$ if and only if $Du(\{ x \}) \not= 0$, then
$\overline{u}$ is continuous in $I \setminus J_u$ and has a jump discontinuity at any
point of $J_u$:
$$\overline{u}(x_{-})  := \lim_{y \uparrow x}\overline{u}(y) = Du(]0,x[), \ \ \ \
\ \overline{u}(x_{+}) := \lim_{y \downarrow x}\overline{u}(y) = Du(]0,x]) \ \ \ \forall
\, x \in J_u,$$ where by simplicity we are assuming that $ I = ]0,L[$. Consequently,
$$\overline{u}(x_{+}) - \overline{u}(x_{-})
 = Du(\{ x \}) \ \ \ \forall \, x \in J_u.$$
Moreover, $\overline{u}$ is differentiable at ${\mathcal L}^1$ a.e. point of $I$, and the
derivative $\overline{u}'$ is the density of $Du$ with respect to ${\mathcal L}^1$, being $\mathcal{L}^d$ the d-dimensional Lebesgue's measure. For
$u \in BV(I)$, the  measure $Du$ decomposes into its absolutely continuous and singular
parts $Du = D^{ac} u + D^s u$. Then $D^{ac} u = \overline{u}' \ {\mathcal L}^1$. Obviously, if
$u \in BV(I)$ then $u \in W^{1,1}(I)$ if and only if $D^su \equiv 0$, and in this case we
have $Du = \overline{u}' \ {\mathcal L}^1$. From now on, when we deal with pointwise
valued $BV$-functions we always shall use the good representative. Hence, in the case $u
\in W^{1,1}(I)$, we shall assume that  $u \in C(\overline{I})$.

Given $\z \in W^{1,1}(I)$ and $u \in BV(I)$, by $\z Du$ we mean the Radon measure
in $I$ defined as
$$
\langle \varphi, \z Du \rangle := \int_{0}^{L} \varphi \z \, Du \ \
\ \ \ \ \forall \, \varphi \in C_c( ]0,L[).
$$
We need the following
integration by parts formula, which can be proved using a suitable regularization of $u
\in BV(I)$ as in the proof of Theorem C.9. in \cite{ACMBook}.

\begin{lemma}\label{IntBP} If $\z \in W^{1,1}(I)$ and $u \in
BV(I)$, then
\begin{equation}\nonumber\label{EIntBP}
\int_{0}^{L} \z Du + \int_{0}^{L} u(x) \z^{\prime}(x) \, dx = \z(L) u({L}_{-})- \z(0)
u({0}_{+}).
\end{equation}
\end{lemma}

\subsection{Properties of the Lagrangian}

 Hereafter $C$ denotes a generic constant, its value may change from line to line.
 We define
\begin{equation}\label{campo1}
\a(z, \xi):= \frac{\nu \vert z \vert \xi}{\sqrt{z^2 + \frac{\nu^2}{c^2}\vert \xi
\vert^2}}.
\end{equation}
 We assume $\a(z,0)=0$ for all $z \in \RR$.
Then $ \a(z, \xi) = \partial_{\xi}F(z, \xi)$, being the Lagrangian
 $$F(z, \xi) := \frac{c^2}{\nu}\vert z \vert \sqrt{z^2 + \frac{\nu^2}{c^2}\xi^2}.$$
  By the convexity of $F$,
 \begin{equation}
 \label{convex_F}
 \a(z,\xi) (\eta - \xi) \le F(z,\eta) - F(z,\xi) \ \ \ \ \hbox{for all} \ z, \xi, \eta \in \R
 \end{equation}
 Note that we have
 \begin{equation}\label{BOUNFG}
c \vert z \vert  \vert \xi \vert - \frac{c^2}{\nu} z^2 \leq \a(z, \xi) \xi \leq c M \vert
\xi \vert \ \ \ \ \hbox{for all} \ z, \xi \in \R, \ \ \vert z \vert \leq M.
 \end{equation}
Moreover, using (\ref{convex_F}) it is easy to see that
\begin{equation}\label{monotoneineq}
(\a(z,\xi)-\a(z, \hat{\xi})) \cdot (\xi - \hat{\xi}) \geq 0
\end{equation}
for any $(z,\xi), (z,\hat{\xi})\in \R\times \R$, $\vert z\vert \leq
M$.

 We introduce the following notation to ease the way in which our functional calculus is
written: for any function $q$ let $J_q(r)$ denote its primite, i.e.,
 $J_q(r)=\int_0^r q(s) \, ds$.

Assume that $f:\R\times \R \to [0, \infty[$ is a continuous
function convex in its last variable such that
\begin{equation}\label{LGRWTHnox}
0 \leq f(z, \xi)  \leq C(1+ | \xi |) \qquad \forall (z,\xi)\in
\R\times\R, \,  |z|\le M.
\end{equation}
for some constant  $C \geq 0$ which may depend on $M$. Given $f(z,\xi)$, we define its
recession function as
$$
 f^0(z, \xi) = \lim_{t \to 0^+} tf \left( z, \frac{\xi}{t} \right).
$$

We assume that $f^0(z,
\xi) = \varphi(z) \psi^0(\xi)$, with $\varphi$ Lipschitz continuous, $\psi^0$ homogeneous of degree $1$.
Then,
working as in \cite{ACMRelatE}, if for a fixed function
$\phi \in C_c(]0,L[)$ we define the operator $\mathcal{R}_{\phi f} : BV(]0,L[) \rightarrow \R$ by
\begin{equation}\label{funct1}
\mathcal{R}_{\phi f}(u):= \int_0^L \phi(x) f(u(x), u^{\prime}(x)) \, dx
 + \int_0^L \phi(x) \psi^0 \left(\frac{D u}{|Du|} \right) \vert D^sJ_\varphi(u) \vert,
\end{equation}
 we have that $\mathcal{R}_{\phi f}$ is lower semi-continuous
respect to the $L^1$-convergence.

For instance,  we discuss here for future usage one of the most
recurrent cases: define $\theta(z)= c|z|$, and note that $F^0(z,
\xi) =  \theta(z) \psi^0(\xi)$, with $\psi^0 (\xi)=\vert \xi
\vert.$ Therefore,
$$
\mathcal{R}_{\phi F}(u):= \int_0^L \phi(x) F(u(x), u^{\prime}(x)) \, dx +
\frac{c}{2}\int_0^L\phi(x) \vert D^s(u^2) \vert
$$
is lower semi-continuous in $BV(]0,L[)$
respect to the $L^1$-convergence.

We shall consider the function $h : \RR \times \RR \rightarrow \RR$ defined by
$
    h(z,\xi):= \a(z,\xi)\cdot  \xi.
$
Note that
\begin{equation}
\label{positivity} h(z,\xi) \ge 0 \ \ \ \ \forall \xi, z \in \RR.
\end{equation}
We will make use of the following property:
\begin{equation}
\label{recession} h^0(z,\xi)=F^0(z,\xi) \quad \forall \xi, z \in \RR.
\end{equation}

As for the Dirichlet problem (see \cite{ACMMRelatDirich}), in
general, the data in $L$ is not taken pointwise; we need to introduce functionals
 that take into account the boundary. The following result is a particular case of Theorem 2.4 in \cite{ACMMRelatDirich}

\begin{theorem}\label{thlsc}
  Let $f$ be verifying (\ref{LGRWTHnox}) and $f^0(z,
\xi) = \varphi(z) \vert \xi \vert$, with $\varphi$ Lipschitz
continuous, let $\phi\in C([0,L])^+$ be given. Then, the
functional
  $\mathcal F_{\phi f}^0:BV(]0,L[)\longrightarrow \R$
  defined by $$\mathcal
  F_{\phi f}^0(u):=\mathcal R_{\phi f}(u)+ \phi(L) \, |J_{\varphi}(u)(L_-)| $$ is lower semi-continuous
  with respect to the $L^1-$convergence.
\end{theorem}

\subsection{Spaces of truncated functions and associated calculus}
We need to consider the following truncature functions. For $a < b$, let $T_{a,b}(r) :=
\max(\min(b,r),a)$. As usual, we denote $T_k = T_{-k, k}$. We also  consider the
truncature functions  $T^l_{a,b}(r):= T_{a,b}(r) - l$  ($l \in \R$). We denote
$$\mathcal T_r:= \{ T_{a,b} \ : \ 0 < a < b \}, \ \ \
{\mathcal T}^+:= \{ T^l_{a,b} \ : \ 0 < a < b, \ l \in \R , \,
T^l_{a,b} \geq 0 \}.$$ Given any truncature function $T_k$, we define
$$T_k(r)^+:=\max\{T_k(r),0\} \ \ \hbox{and}  \ \ \ T_k(r)^-:=\min\{T_k(r),0\}=-T_k(-r)^+, \
 r\in\R.$$

Consider the function space $$TBV^+(I):= \left\{ u \in L^1(I)^+  \ :  \ \ T(u) \in BV(I),
\ \ \forall \ T \in \mathcal T_r \right\};$$ we want to give a sense to the Radon-Nikodym
derivative $ u'$ of a function $u \in TBV^+(I)$. Using chain's rule for BV-functions
(see, for instance, \cite{Ambrosio}), and with a similar proof to the one given in Lemma
2.1 of \cite{Benilanetal}, we obtain the following result.

\begin{lemma}\label{WRN}
For every $u \in TBV^+(I)$ there exists a unique measurable function $v : I \rightarrow
\R$ such that
\begin{equation}\label{E1WRN}
(T_{a,b}(u))' = v \1_{[a < u  < b]} \ \ \ \ \ {\mathcal L}^1-{\rm a.e.}, \ \ \forall \
T_{a,b} \in \mathcal T_r.
\end{equation}
\end{lemma}

Thanks to this result we define $u'$ for a function $u \in
TBV^+(I)$ as the unique function $v$ which satisfies
(\ref{E1WRN}). This notation will be used throughout in the
sequel. The notation $\partial_x$ will also be used in the case of
functions of several
 variables (say $t$ and $x$), for the same purposes, whenever there is some risk of
 confusion.

 We denote by ${\mathcal P}$ the set of Lipschitz continuous function $p : [0, +\infty[ \rightarrow \R$
satisfying $p^{\prime}(s) = 0$ for $s$ large enough, and write ${\mathcal P}^+:= \{ p \in
{\mathcal P} \ : \ p \geq 0 \}$. We recall the following result (\cite{ACM4:01}, Lemma
2).

\begin{lemma}
\label{calc_producto} If $u \in TBV^+(I)$, then $p(u) \in BV(I)$
for every $p \in {\mathcal P}$ such that there exists $a > 0$ with
$p(r) = 0$ for all $0 \leq r \leq a$. Moreover, with the above
notation $[p(u)]'=p'(u) u' \ \ \mathcal{L}^1$-a.e.
\end{lemma}

For $u \in TBV^+(]0,L[)$ we will define
$$u(0_+):= \lim_{n\to \infty} T_{\frac{1}{n},n}(u)(0_+) \ \ \ \hbox{and} \ \ \ u(L_-):= \lim_{n\to \infty} T_{\frac{1}{n},n}(u)(L_-).$$It is
easy to see that the above limits exist.

 Let $S \in \mathcal{P}^+$ and $T = T_{a,b}^a$. Given $u \in TBV^+(]0,L[)$,
 Lemma \ref{calc_producto} assures that
$
   S(u)T(u), J_{T'S}(u),J_{TS'}(u) \in BV(]0,L[).
$
Moreover,
$
   D(S(u)T(u))=DJ_{T'S}(u) + DJ_{TS'} (u)
$
and hence, if $\z \in W^{1,1}(]0,L[)$,
$$
   \z D(T(u)S(u)) = \z DJ_{T'S}(u) + \z DJ_{TS'}(u).
$$

For $u \in TBV^+(]0,L[)$, $\phi \in C_c(]0,L[)$, $T=T_{a,b} - l
\in \mathcal{T}^+$ and $f$ as in the previous subsection -see
(\ref{LGRWTHnox})-, we define the functional
$$
  \mathcal{R}(\phi f,T)(u):= \mathcal{R}_{\phi f}(T_{a,b}(u)) + \int_{[u \le a]} \phi(x) (f(u(x),0)-f(a,0)) \, dx
$$
$$
 - \int_{[u \ge b]} \phi(x) (f(u(x),0)-f(b,0)) \, dx.
$$
We have that $\mathcal{R}(\phi f,T)(\cdot)$ is lower-semi-continuous in $TBV^+(]0,L[)$
 with respect to the $L^1$-convergence.

Given $S, T \in \mathcal{T}^+$ and $u \in TBV^+(]0,L[)$, we define
the following Radon measures in $]0,L[$,
$$
  \langle F(u, DT(u)), \phi \rangle: = \mathcal{R}(\phi F,T)(u),
\quad
 \langle F_S(u, DT(u)), \phi \rangle : = \mathcal{R}(\phi S F,T)(u),
$$
$$
\langle  h(u, DT(u)), \phi \rangle: = \mathcal{R}(\phi h,T)(u), \quad
 \langle  h_S(u, DT(u)), \phi \rangle : = \mathcal{R}(\phi S h,T)(u),
$$
for $\phi \in C_c(]0,L[)$.
Using (\ref{funct1}) and (\ref{recession}), we compute
 $$
    F(u,DT(u))^s = \frac{c}{2}\left|D^s(T(u))^2 \right| = h(u,DT(u))^s, \ \ \ F_S(u,DT(u))^s = \left|D^sJ_{S \theta}(T(u)) \right| = h_S(u,DT(u))^s,
 $$
 $$
        h(u,DT(u))^{ac}= h(u,(T(u))'), \ \ \  h_S(u,DT(u))^{ac}= S(u) h(u,(T(u))').
 $$

\section{The Elliptic Problem}

Given $v \in L^1(]0,L[)$, we are interested in the following
problem:
\begin{equation}\label{Elliptictproblem}\left\{
\begin{array}{ll}
 - \left(\a(u,u') \right)'  =  v & \ \hbox{in} \ \  ]0, L[
\\ \\ -  \a(u,u')|_{x=0}  = \beta > 0 \ \ \hbox{and} \ \ u(L)=0,
\end{array}
\right.
\end{equation}
where $\a$ is given by (\ref{campo1}). We introduce the following concept of solution for
problem (\ref{Elliptictproblem}).

\begin{definition}\label{DefSPE} {\rm Given $v \in L^1(]0,L[)$,
we say that $u \geq 0$ is an {\it entropy solution} of
(\ref{Elliptictproblem}) if $u \in TBV^+(]0,L[)$ and $\a(u,
 u') \in C([0,L])$ both satisfy
 \begin{equation}\nonumber\label{Cond1}
 v  = - D \a(u, u') \ \ \ \ \ \ {\rm in} \ \ {\mathcal
D}^{\prime}(]0,L[),
\end{equation}
\begin{equation}\nonumber\label{Cond2}
 - \a(u, u')(0) =  \beta, \ \ \ \ \ {\rm and} \ \ \ \
  \a(u, u')(L) = - c u(L_{-}).
\end{equation}
\begin{equation}\label{Cond2,2} h(u,DT(u))\le \a(u,u')DT(u)
\ \ \ {\rm as \ measures} \ \ \forall \, T \in {\mathcal T}^+
\end{equation}
\begin{equation}\label{SCond2}
\displaystyle   h_S(u,DT(u))\le \a(u,u')DJ_{T'S}(u) \ \ \ \hbox{as  measures}
\\ \\  \ \displaystyle\forall \, S \in {\mathcal P}^+,  \, T\in {\mathcal T}^+.
\end{equation}
}

\end{definition}

\noindent
Note that (\ref{Cond2,2}) can be rewritten as $h(u,DT(u))^s \le [\a(u,u') DT(u)]^s$, and
thus it is equivalent to
\begin{equation}\nonumber\label{Cond2,2Bis}
\frac{c}{2} \vert D^s((T(u))^2) \vert \leq \a(u, u') D^sT(u) \ \ \ {\rm as \ measures} \
\ \forall \, T \in {\mathcal T}^+.
\end{equation}
Also we have that (\ref{SCond2}) can be rewritten as
$h_S(u,DT(u))^s \le [\a(u,u') DJ_{T'S}(u)]^s$, and is equivalent
to
\begin{equation}\nonumber\label{SCond2Bis}
\displaystyle \vert D^s(J_{S \theta}(T(u)) \vert \leq \a(u, u') D^sJ_{T'S}(u) \ \ {\rm as
\ measures} \
 \ \ \forall \, S\in {\mathcal P}^+,  \, T\in {\mathcal T}^+.
 \end{equation}
Observe that since $- \a(u,
 u')(0) =  \beta$, we have
 \begin{equation}\label{uencero}
u(0_+) \geq \frac{\beta}{c} > 0.
 \end{equation}

We introduce now the main result of this section.
\begin{theorem}\label{EUTElliptic}
 For any
 $0 \leq f \in L^{\infty}(]0,L[)$
 there exists a unique
 entropy solution $u \in TBV^+(]0,L[)$ of the
 problem
\begin{equation}
\label{Elliptictproblem1}
\left\{
\begin{array}{ll}
u - \left(\a(u,u') \right)'  =  f & \hbox{in} \  ]0, L[
\\ \\ -  \a(u,u')|_{x=0}  = \beta > 0 \quad & u(L)=0,
\end{array}
\right.
\end{equation}
which satisfies $\|u\|_\infty \le  M(\beta, c, \nu, \|f\|_\infty)$.

\noindent
Moreover, let $u, \overline{u}$ be two  entropy solutions of (\ref{Elliptictproblem1}) associated to $f, \overline{f} \in L^{1}(]0,L[)^+$,
respectively. Then,

$$\int_0^L (u - \overline{u})^+ \, dx \leq \int_0^L (f - \overline{f})^+\, dx.$$
\end{theorem}
\noindent{\bf Proof.} {\it Existence of entropy solutions.} We
divide the existence proof in different steps.

\noindent{\it Step 1. Approximation and basic estimates.} Let $0
\leq f \in L^{\infty}(]0,L[)$. For every $n \in \NN$, consider
$\a_n(z , \xi)$ $:= \a (z, \xi) + \displaystyle \frac{1}{n} \xi$. As
a consequence of the results about pseudo-monotone operators
in \cite{Browder} we know that $\forall n \in \N$ there exists a
unique $ u_n \in W^{1,2}(]0,L[)$ such that $u_n(L)=0$ and
\begin{equation}
   \int_0^L u_n v \, dx +  \int_0^L \a(u_n,u_n^{\prime}) v^{\prime} \, dx + \frac{1}{n}\int_0^L u_n^{\prime} v^{\prime} \, dx -
   \beta v(0) = \int_0^L f v \, dx
\label{weak form}
\end{equation}
for all $v \in W^{1,2}(]0,L[), \, v(L)=0$.

 The following result can be easily obtained by  multiplication by $u_n^-$ and
 integration over $[0,L]$.

\begin{lemma}\label{f1}
 The functions $u_n$ are non-negative  $ \forall \, n \in \NN$.
\end{lemma}

 Now we give a bound for the sequence $u_n$ at zero.

\begin{lemma}
\label{Acotcero}
The sequence $\{u_n(0) \}$ is bounded. More precisely,
\begin{equation}
\nonumber
0 \leq u_n(0) \leq \left\{ \begin{array}{l} \displaystyle
\frac{4\beta c}{\nu} + \sqrt{\frac{2cL}{\nu}} \Vert f
\Vert_{\infty} \quad
\quad \hbox{if} \ \ c > \sqrt{\nu} \\ \\
\displaystyle \frac{4 \beta}{c} + \sqrt{\frac{2L}{c}} \Vert f
\Vert_{\infty}  \quad \quad \quad  \hbox{if} \ \ c \leq
\sqrt{\nu}.
\end{array} \right.
\end{equation}
\end{lemma}
\noindent{\bf Proof.}
Taking  $v=u_n$ in (\ref{weak form}), we get
\begin{equation}
\label{EE33}
\int_0^L \left(u_n^2 + \a(u_n,u_n^{\prime}) u_n^{\prime}
+ \frac{1}{n} ((u_n)^{\prime})^2 \right) \, dx = \beta u_n(0) + \int_0^L
f u_n \, dx.
\end{equation}
Then, dropping non-negative terms
and using Young's inequality, we get
\begin{equation}\label{EE30o3}
 \int_0^L u_n^2 \, dx \leq \int_0^L f^2\, dx + 2 \beta \, u_n(0).
\end{equation}
Now we can write
$u_n|u_n^{\prime}|=\frac{1}{2}|(u_n^2)^{\prime}|$, and taking
into account (\ref{BOUNFG}) we have  $u_n^{\prime}
\a(u_n,u_n^{\prime}) \ge c u_n|u_n^{\prime}|- \frac{c^2}{\nu}u^2_n$. Then, from (\ref{EE33}), we obtain
\begin{eqnarray}
\label{Ebuena0}
 \int_0^L  \left( \frac{c}{2} |(u_n^2)^{\prime}| + \frac{((u_n)^{\prime})^2}{n}  \right)
dx \leq \int_0^L \left(  \left(\frac{c^2}{\nu} -1 \right) u_n^2 +
f u_n \right)  dx + \beta u_n(0).
 \end{eqnarray}

Assuming now that $\frac{c^2}{\nu} -1 > 0$, we apply Young's inequality in the right hand side of
(\ref{Ebuena0}), which now reads
\begin{equation}\nonumber
\label{Ebuena1}
\left(\frac{c^2}{\nu} -\frac{1}{2} \right) \int_0^L u_n^2 \, dx + \frac{1}{2} \int_0^L f^2 \, dx +
\beta u_n(0).
 \end{equation}
As $c > \sqrt{\nu}$ we have $\frac{c^2}{\nu} -\frac{1}{2}>0$, which allows us to bring in (\ref{EE30o3}), thus obtaining
\begin{equation}\label{Ebuena}
    \frac{c}{2} \int_0^L |(u_n^2)^{\prime}| \, dx + \frac{1}{n}\int_0^L ((u_n)^{\prime})^2\, dx
    \leq \frac{c^2}{\nu}\int_0^L f^2 \, dx+  \frac{2\beta c^2}{\nu} u_n(0).
    \end{equation}
Then, we have
$$
  \frac{c}{2}|u_n^2(0)| = \frac{c}{2}|u_n^2(L)-u_n^2(0)| =
  \frac{c}{2} \left\vert \int_0^L (u_n^2)^{\prime} \, dx
  \right\vert \leq \frac{c^2}{\nu}\int_0^L f^2 \, dx+ \frac{ 2\beta c^2}{\nu} u_n(0),
$$
from where we get that $u_n^2(0) - \frac{4 \beta c}{\nu} u_n(0) - \frac{2c}{\nu} \Vert f \Vert_2^2 \leq 0.$
Hence, for all $n \in \NN$,
\begin{equation}\nonumber
\label{Ebuena44}
0 \leq u_n(0) \leq \frac{1}{2} \left(\frac{4 \beta c}{\nu} +\sqrt{\left(\frac{4 \beta c}{\nu}\right)^2 + \frac{8c}{\nu}
 \Vert f \Vert_2^2   }
\right) \leq \frac{4\beta c}{\nu} + \sqrt{\frac{2c}{\nu}} \Vert f
\Vert_2.
\end{equation}

In case that $c^2/\nu -1 \le 0$, from
(\ref{Ebuena0}) we obtain
$$
\frac{c}{2} \int_0^L |(u_n^2)^{\prime}|\, dx + \frac{1}{n}\int_0^L
((u_n)^{\prime})^2\, dx \leq  \int_0^L f u_n \, dx + \beta u_n(0).
$$
Then, using Young's inequality and having in mind (\ref{EE30o3}),
we get
\begin{equation}
\label{EbuenaQQ}
    \frac{c}{2} \int_0^L |(u_n^2)^{\prime}| \, dx + \frac{1}{n}\int_0^L ((u_n)^{\prime})^2\, dx
\leq  \int_0^L f^2 \, dx+  2\beta u_n(0).
    \end{equation}
Thus, we have that $u_n^2(0)-  \frac{4}{c} \beta u_n(0) - \frac{2}{c} \int_0^L f^2 \leq 0$, from where it follows that for all $n \in \NN$,
\begin{equation}\nonumber
\label{Ebuena55}
0 \leq u_n(0) \leq \frac{1}{2}
\left(\frac{4 \beta}{c} + \sqrt{\left(\frac{4 \beta}{c}\right)^2 +
\frac{8}{c}
 \Vert f \Vert_2^2   }
\right) \leq \frac{4 \beta}{c} + \sqrt{\frac{2}{c}} \Vert f
\Vert_2.
\end{equation}

\hfill$\Box$

By (\ref{Ebuena}), (\ref{EbuenaQQ}) and Lemma \ref{Acotcero}, we get
\begin{equation}\label{Ebuenisi}
 \frac{c}{2} \int_0^L |(u_n^2)^{\prime}| \, dx +
    \frac{1}{n}\int_0^L ((u_n)^{\prime})^2 \, dx \leq C \ \ \ \ \forall \, n \in \NN.
\end{equation}

\begin{lemma}
\label{BOUNDEDD}
The sequence $ \{u_n \ : \ n\in \NN \}$  is uniformly bounded in $  L^{\infty}(0,L)$.
\end{lemma}
\noindent{\bf Proof.} By Lemma \ref{Acotcero}, we know that $
M=\mbox{max} \left\{\|f\|_\infty, \mbox{max}\{u_n(0) \ : \ n \in
\NN\} \right\}$ is finite. Then, taking $v = \left(u_n -  M
\right)^+$ as test function in (\ref{weak form}),
 it is easy to see that $
\Vert u_n \Vert_{\infty} \leq M $
and Lemma \ref{BOUNDEDD} holds.
\hfill$\Box$

\begin{lemma}
The sequence $\{ u_n\}$ is uniformly bounded in $TBV^+(]0,L[)$. Furthermore, there exists a function $0\le u \in TBV^+(]0,L[) \cap L^\infty (]0,L[)$ such that (up to subsequence) $u_n \to u$ a.e. and strongly in $L^1(]0,L[)$.
\end{lemma}
\noindent{\bf Proof.}
By Lemma \ref{BOUNDEDD}, extracting a subsequence if necessary,
we may assume
 that $u_{n}$ converges weakly in $L^2(]0,L[)$
 to some non-negative function $u$ as $n \rightarrow + \infty$.
Moreover, by Lemma \ref{BOUNDEDD} again, we have that  $0 \leq u \in
L^{\infty}(]0,L[)$. On the other hand, if $0 < a < b$, by the
coarea formula and (\ref{Ebuenisi}), we have
$$
   \int_0^L | (T_{a,b} (u_n))^\prime| \, dx = \int_a^b \big|D \1_{[u_n\le t]}\big|(]0,L[) \, dt
   = \int_a^b \big|D \1_{[u_n^2 \le t^2]}\big|(]0,L[) \, dt
$$
$$
 =  \int_{a^2}^{b^2} \big| D \1_{[u_n^2 \le s]}\big|(]0,L[) \frac{ds}{2\sqrt{s}} \le
   \frac{1}{2a}\int_0^L |(u_n^2)^{\prime}| \, dx \le \frac{C}{ a}.
$$
 Consequently, we may assume that $u_{n}$ converges  almost everywhere to $u$.
Then, by the Vitali Convergence Theorem, we get that $u_n \to u$ in $L^1(]0,L[)$, and
using the above estimate on the gradients we obtain that $u \in TBV^+(]0,L[)$.
\hfill$\Box$

Since $|\a(u_n,u^{\prime}_n)| \le c |u_n|$, by Lemma \ref{BOUNDEDD} we
may assume that
\begin{equation}
\label{weakz3}
\a(u_n, u^{\prime}_n)  \rightharpoonup \z \quad
\hbox{as $n \to \infty$, weakly$^*$ in $L^{\infty}(]0,L[)$}.
\end{equation}
By assumption we have that $\a(u_n,u^{\prime}_n) = c \vert
u_n \vert \b(u_n,u^{\prime}_n)$ with $\vert \b(u_n,u^{\prime}_n)
\vert \leq 1$ (independent of $n$), $\Vert u_n\Vert_\infty \leq
 M$, and $u_n \to u$ a.e. as $n\to \infty$, so we may assume that $
\b(u_n, u^{\prime}_n)  \rightharpoonup \z_b$ as $n \to
\infty$, weakly$^*$ in $L^\infty(]0,L[)$,
and
\begin{equation}\label{bbCONL2}
\z = cu \z_b, \ \ \ \ \ \ \  {\rm with} \ \ \ \ \Vert \z_b \Vert_\infty  \leq 1.
\end{equation}
On the other hand,  by (\ref{Ebuenisi}),
\begin{equation}\label{CONL2}
\frac{1}{n}  u^{\prime}_n \to 0 \ \ \ \ {\rm in} \ \ L^2(]0,L[).
\end{equation}

Given $\phi \in {\mathcal D}(]0,L[)$, taking $v = \phi$ in
(\ref{weak form}) we obtain
$$
   \int_0^L u_n \phi \, dx +  \int_0^L \a(u_n,u_n^{\prime}) \phi^{\prime} \, dx +
   \frac{1}{n}\int_0^L u_n^{\prime} \phi^{\prime}\, dx = \int_0^L f \phi \, dx
$$
Letting $n \rightarrow + \infty$, having in mind (\ref{weakz3})
and (\ref{CONL2}), we obtain
\[
\int_0^L (f - u) \phi \ dx = \int_0^L \z \cdot  \phi^{\prime} \
dx,
\]
that is,
\begin{equation}\label{div}
f - u = - D\z, \ \ \ \ \ \ {\rm in} \ \ {\mathcal
D}^{\prime}(]0,L[)
\end{equation}
and
\begin{equation}\nonumber\label{ConvDiv}
(\a_n(u_n,  u_n^{\prime}))^{\prime} \rightharpoonup D\z \ \ \ \
{\rm weakly \ \ in} \ \ L^2(]0,L[).
\end{equation}
Note that by (\ref{div}), we have $\z \in W^{1,1}(]0,L[)$ and $D\z = \z^{\prime}$.

Working as in the proof of Lemma 4.2 of \cite{ACMRelatE}, we can prove the identification
 \begin{equation}
 \z(x)=\a(u(x),u^\prime(x)) \ \ \mbox{a.e.} \, \, x \in ]0,L[
 \label{ident_campo}
 \end{equation}
From (\ref{ident_campo}) and (\ref{div}) it follows that
\begin{equation}\nonumber
f- u = - D \a(u,u^\prime), \quad \mbox{in}\, \mathcal{D}'(]0,L[)
\label{limit_formulation}
\end{equation}

\begin{lemma}
The flux $-\a (u,u')$ verifies the Neumann condition at $x=0$.
\end{lemma}
\noindent{\bf Proof.}
    Let $w \in W^{1,1}(]0,L[)$ such that $w(L)=0$ and consider $w_k \in W^{1,2}(]0,L[)$ with  $w_k(L) =0$ for all
     $k \in \NN$, $w_k \to \hat{w}$  pointwise and $w_k^{\prime} \to  w^{\prime}$
     in $ L^1(]0,L[)$. Taking in (\ref{weak form}) $w_k$ as test function and letting $n \to + \infty$,
    we get
    $$\int_0^L u w_k \, dx +  \int_0^L \z w_k^{\prime} \, dx -
   \beta w_k(0) = \int_0^L f w_k \, dx.$$
   Then, letting $k \to +\infty$ we arrive to
    \begin{equation}
   \int_0^L u w \, dx +  \int_0^L \z w^\prime \, dx -
   \beta w(0) = \int_0^L f w \, dx.
\label{Cedsr}
\end{equation}

Fixed $w \in BV(]0, L[)$ such that $w(L_{-}) = 0$, let $w_m \in
W^{1,1}(]0, L[)$ with $w_m(L) = 0$, $w_m (0) = w(0_+)$, and such
that $w_m \to w$ in $ L^1(]0, L[).$
Taking in (\ref{Cedsr}) $w_m$ as test functions and integrating by
parts we get
\begin{eqnarray*}
    \int_0^L (f-u) w_m \, dx =   \int_0^L \z w_m^{\prime} \, dx -
   \beta w(0_+) = -  \int_0^L \z^{\prime} w_m \, dx - w(0_+) (\z(0)+
   \beta ),
\end{eqnarray*}

and letting $m \to + \infty$, we obtain
$-\z(0) = \beta.$
\hfill$\Box$

\begin{lemma}
Let $S \in \mathcal{P}^+$, $T\in \mathcal T^+$ and $\phi\in C^1([0,L])$, $\phi\geq 0$, with $\phi(0) = 0$. Then
\begin{equation}
\label{Interm2}
\begin{array}{l}
  \displaystyle
\int_0^L \phi F(u,DT(u)) + \phi(L) \frac{c}{2} \vert
(T(u))^2(L_{-}) \vert
\\ \\
\leq \displaystyle \int_0^L \phi \z DT(u) +
\int_0^L \phi
  F(u,0) \, dx - \phi(L) T(u)(L_{-})  + \phi(L) \left| J_{\theta }(T(0))\right|
  \end{array}
\end{equation}
and
\begin{equation}
\label{S2}
\begin{array}{l}
 \displaystyle
\int_0^L \phi F_{S}(u,DT(u))+\phi(L) \, |J_{\theta S}(T(u)(L_{-}))|
 \\ \\
 \leq
\displaystyle\int_0^L \phi \z DJ_{T'S}(u) +\int_0^L \phi S(u)
F(u,0) \, dx
\\ \\
\displaystyle - \phi(L) \z(L)J_{T'S}(u(L_-))
 + \phi(L) \left| J_{\theta S}(T(0))\right|.
\end{array}
\end{equation}
In particular,

\begin{equation}\label{Interm2*}
  F(u,DT(u))\leq \z DT(u)+
  F(u,0)\mathcal
  L^{1} \ \ \ \hbox{\rm as measures in }]0,L[.
\end{equation}
\begin{equation}
\label{S2*}
 \displaystyle F_S(u,DT(u)) \displaystyle \leq \z D(J_{T'S}(u)))+
S(u) F(u,0) \, \mathcal{L}^1  \ \hbox{\rm as measures in }]0,L[.
\end{equation}

\end{lemma}
\noindent{\bf Proof.}
 We will only prove (\ref{S2}), the proof of (\ref{Interm2}) being similar. Let $0 \leq \phi \in {C^1}([0,L])$ with $\phi(0) = 0$.

Since $\mathcal{F}_{\phi SF}^0$ is l.s.c. with respect to the
$L^1$-convergence, letting $n\to\infty$ we obtain
\[\begin{split}
& \int_0^L \phi F_{S}(u,DT(u))+\phi(L) \, |J_{\theta S}(T(u)(L_{-}))|
\\
& \leq\liminf_{n\to\infty}\int_0^L \phi S(u_n)F(u_n, T(u_n)^{\prime})\
dx + \phi(L) \left| J_{\theta S}(T(0))\right|
\\
& \leq\limsup_{n\to\infty}\int_0^L \phi S(u_n)F(u_n,T(u_n)^{\prime})\ dx + \phi(L) \left| J_{\theta S}(T(0))\right|
\end{split}\]
 By the convexity (\ref{convex_F}) of $F$ and using that
$
\a(u_n, T(u_n)^{\prime}) T(u_n)^{\prime} =  \a(u_n,u_n^{\prime}) T(u_n)^{\prime},
$
we have
\[\begin{split}
& \int_0^L \phi S(u_n) F(u_n, T(u_n)^{\prime})  \, dx
\\
 & \leq \int_0^L
\phi S(u_n) \a(u_n, T(u_n)^{\prime}) T(u_n)^{\prime}  \ dx  +
\int_0^L \phi  S(u_n) F(u_n,0)
 dx
\\
& = \int_0^L \phi \a(u_n, u_n^{\prime})
(J_{T'S}(u_n))^{\prime}  \, dx  +  \int_0^L \phi  S(u_n) F(u_n,0)
 dx.
\end{split}\]
Now we take $v =
J_{T'S}(u_n)\phi$ as test function in (\ref{weak form})) and
  we obtain
\[\begin{split}
 & \int_0^L \phi \a(u_n, u_n^{\prime}) (J_{T'S}(u_n))^{\prime}
 \, dx  + \frac{1}{n} \int_0^L \phi u_n^{\prime}
(J_{T'S}(u_n))^{\prime} \, dx
\\
& =  \int_0^L (f -
u_n)J_{T'S}(u_n) \phi \, dx - \int_0^L J_{T'S}(u_n) \a(u_n,
u_n^{\prime}) \phi^{\prime} \, dx - \frac{1}{n} \int_0^L
J_{T'S}(u_n) u_n^{\prime} \phi^{\prime} \, dx.
\end{split}\]
Letting $n\to\infty$ we get
\[\begin{split}
\limsup_n \int_0^L \phi \a(u_n,u_n^{\prime})
(J_{T'S}(u_n))^{\prime}
 \, dx  \leq & \int_0^L \phi (f-u)J_{T'S}(u) \, dx -
\int_0^L J_{T'S}(u) \z  \phi^{\prime} \, dx
\\
= & \int_0^L \phi \z D(J_{T'S}(u)) - \phi(L) \z(L)J_{T'S}(u(L_-)).
\end{split}\]
Finally,
\[\begin{split}
& \int_0^L \phi F_{S}(u,DT(u))+\phi(L) \, |J_{\theta S}(T(u))(L_{-})| \leq\int_0^L
\phi \z DJ_{T'S}(u)
\\
& + \phi(L) \left| J_{\theta S}(T(0))\right| - \phi(L) \z(L)J_{T'S}(u(L_-))+\int_0^L \phi S(u) F(u,0) \, dx
\end{split}\]
and
(\ref{S2}) holds.

\hfill$\Box$

\begin{lemma}
The inequalities  (\ref{Cond2,2}) and (\ref{SCond2})
hold.
\end{lemma}

\smallskip\noindent {\bf Proof.} Using (\ref{Interm2*}) and the fact that $h(u,DT(u))$ is a measure
concentrated in $]0,L[$, it follows that
 $$h(u, DT(u))^s= F(u, DT(u))^s  \leq  (\z \
 DT(u))^s.$$
Hence,
$$ \z DT(u)) = \z  T(u)^{\prime}\mathcal L^1 + (\z
 DT(u))^s \geq \z T(u)^{\prime}\mathcal L^1 + h(u, DT(u))^s = h(u, DT(u)),$$
and  (\ref{Cond2,2}) holds.

\medskip

Using (\ref{S2*}) we have
$$
\z D(J_{T'S}(u)))  = (\z D(J_{T'S}(u)))^{ac} + (\z
D(J_{T'S}(u)))^s  \geq \z (J_{T'S}(u))^{\prime} + (F_S(u,
DT(u)))^s
$$$$=\z (J_{T'S}(u))^{\prime}\mathcal L^N + (h_S(u,DT(u)))^s =
h_S(u,DT(u)),$$ and we obtain (\ref{SCond2}).\hfill$\Box$

\begin{lemma}\label{lemaBCCD}
The Dirichlet condition $\a(u, u')(L) =-  c u(L_{-})$ holds.
\end{lemma}
\noindent {\bf Proof.}
 Firstly, observe that by (\ref{bbCONL2}) we have
 $$
 \vert \z(L) \vert \leq c u(L_{-}).
 $$
 Then, it is enough to prove the lemma in the case $u(L_{-}) > 0$.
 In that case, again by (\ref{bbCONL2}) and having in mind that $\z$ is continuous in $[0,L]$, we have
\begin{equation}
\label{bvcCONL2}
\z(L) = cu(L_{-}) \xi, \ \ \ \ \ \ \ {\rm with} \ \ \ \vert \xi  \vert \leq 1.
\end{equation}
 Given $T\in T^+$, for $m>1$ we consider
 $S:=T^{m-1}\in \mathcal{P}^+$. Taking singular parts in (\ref{S2}) we have
\begin{equation}
\label{utilx0}
|J_{\theta T^{m-1}}(T(u))(L_{-})|\leq -\z(L) J_{T^{m-1}T'}(u(L_{-})) + \left| J_{\theta T^{m-1}}(T(0))\right|.
\end{equation}
Consider now $T=T_{d,d'}$ with $0 < d \le u(L_-)\le \|u-\|_\infty \le
d'$. Using (\ref{bvcCONL2}), the inequality (\ref{utilx0}) particularizes to
$$
  \frac{c}{2} d^{m+1} +\frac{c}{m+1} \left(u^{m+1}(L_-) - d^{m+1} \right) \le \frac{c}{2} d^{m+1} - \frac{c}{m} \xi u(L_-) \left(u^m(L_-) - d^m \right)
$$
and letting  $d\to 0^+$ we have
$$
\frac{c}{m+1}u^{m+1}(L_{-})  \leq - \frac{c}{m} u((L_{-})) \xi u^m((L_{-}).
$$
Then, since $u(L_{-}) > 0$, we get $\frac{m}{m+1} \leq - \xi$ for
all $1<m$. Therefore, since $\vert \xi  \vert \leq 1$, we have
$\xi = - 1$. Consequently, by (\ref{bvcCONL2}) we finish the
proof.
\hfill$\Box$

\medskip

\noindent{\it Proof of uniqueness.} Let  $u,
\overline{u}$  be  entropy solutions of (\ref{Elliptictproblem1}) associated with  $f, \overline{f} \in L^{1}(]0,L[)^+$, respectively.

Let $\rho_n$ be a classical mollifier in $]0,L[$, $\psi\in
{\mathcal D}(]0,L[)$ and $b
> a
> 2 \epsilon > 0$. Let us write
$$
\xi_{n}(x,y) =  \rho_n(x-y)\psi\left(\frac{x+y}{2}\right), \ \ {\rm
and} \ \ T= T^a_{a,b}.
$$
We need to consider truncature functions of the form $S_{\epsilon, l}(r):= T_{\epsilon}(r - l)^+ = T_{l, l + \epsilon}(r) - l \in {\mathcal T}^+$
 and $S_\epsilon^l(r):= T_{\epsilon}(r - l)^- + \epsilon =
T_{l- \epsilon, l}(r) + \epsilon - l \in {\mathcal T}^+,$ where
$l \geq 0$. Observe that $
S_\epsilon^l(r) = - T_\epsilon(l-r)^+ + \epsilon.
$

 If we denote $\z(y) = \a(u(y),\partial_y
u(y))$ and $\overline{\z}(x) = \a(\overline{u}(x), \partial_x
\overline{u}(x))$, we have
\begin{equation}\nonumber
\label{EE11}
u - \z' = f \ \ \ {\rm and} \ \ \ \overline{u} - \overline{\z}' =
\overline{f} \ \ \ {\rm in} \ \ \ \  {\mathcal
D}^{\prime}(]0,L[).
\end{equation}
Then, multiplying the equation for $u$ by
$T(u(y))S_{\epsilon, \overline{u}(x)}(u(y)) \xi_n(x,y)$, that for
$\overline{u}$ by $T(\overline{u}(x)) S_\epsilon^{u(y)}(\overline{u}(x))
\xi_n(x,y)$, integrating in both variables

\begin{equation}
\label{E1ENTROP}
\begin{array}{l}
\displaystyle\int_0^L \int_0^L
[u(y)T(u(y))-\overline{u}(x)T(\overline{u}(x))]
T_\epsilon(u(y)-\overline{u}(x))^+\xi_n(x,y) \, dxdy
\\ \\
+ \epsilon \displaystyle\int_0^L \int_0^L (\overline{u}(x)-
\overline{f}(x))T(\overline{u}(x)) +
 \xi_n(x,y) \,  dxdy
\\ \\
+ \displaystyle\int_0^L \int_0^L \xi_{n}(x,y) \left( \z
D_y[T(u)S_{\epsilon,\overline{u}(x)}(u)]dx + \overline{\z} D_x[T(\overline{u})
S_\epsilon^{u(y)}(\overline{u})]dy  \right)  \\ \\
+ \displaystyle\int_0^L \int_0^L T(u(y))
S_{\epsilon,\overline{u}(x)}(u(y)) \z(y) \cdot \partial_y \xi_n(x,y)
\, dxdy
\\ \\ + \displaystyle\int_0^L \int_0^L T(\overline{u}(x))
S_\epsilon^{u(y)}(\overline{u}(x)) \overline{\z}(x) \cdot \partial_x
\xi_n(x,y) \, dxdy
\\ \\ = \displaystyle\int_0^L \int_0^L
[f(y)T(u(y))-\overline{f}(x)T(\overline{u}(x))]
T_\epsilon(u(y)-\overline{u}(x))^+\xi_n \, dxdy.
\end{array}
\end{equation}
 Let $I$ denote all the terms at the left hand side of the above identity,
but the first one. From now on, since $u,\z$ are always functions
of $y$, and $\overline{u},\overline{\z}$ are always functions of
$x$, to make our expressions shorter, we shall omit the arguments
except in some cases
where we find useful to remind them.

 With slight modifications of the method used in the proof of uniqueness in \cite{ACMRelatE} we can obtain the following result.
\begin{lemma}
\label{auxiliarstep}
The following inequality is satisfied
\begin{eqnarray*}
\frac{1}{\epsilon}  I \geq o(\epsilon)  -
\displaystyle\int_0^L\left( \int_0^L \xi_n \overline{\z}
D_xT(\overline{u}) \right) dy + \frac{1}{\epsilon}
\int_0^L\int_0^L T_\epsilon (u - \ubar)^+ (T(u) \z - T (\ubar)
\zbar)\cdot (\partial_x \xi_n + \partial_y \xi_n) \, dx \ dy,
\end{eqnarray*}
where $o(\epsilon)$ denotes an expression such that
$o(\epsilon)\to 0$ as $\epsilon\to 0$.
\end{lemma}

By the above lemma, dividing (\ref{E1ENTROP}) by $\epsilon$ and letting
$\epsilon\to 0$ we obtain
$$\int_0^L\int_0^L \xi_n(x,y)(u(y)T(u(y))-\overline{u}(x)T(\overline{u}(x))) {\rm
sign}_0^+(u(y)-\overline{u}(x))\, dx\, dy
$$
$$
+\int_0^L\int_0^L
\rho_n (x-y){\rm sign}_0^+(u(y)-\overline{u}(x))(T(u(y))\z(y)-T(\overline
u(x))\zbar(x))\cdot \psi'\left(\frac{x+y}{2}\right) \, dx \, dy
$$
$$
\leq \int_0^L\int_0^L
\xi_n(x,y)(f(y)T(u(y))-\overline{f}(x)T(\overline{u}(x))) {\rm
sign}_0^+(u(y)-\overline{u}(x))\, dx \,dy + \int_0^L\left( \int_0^L
\xi_n(x,y) \overline{\z} D_xT(\overline{u}) \right) dy,
$$
where
$$
\displaystyle {\rm sign}_0^+(r)=\left\{\begin{array}
  {cc} 1 & \quad \mbox{if } r>0 \\ 0 &\quad \mbox{if } r\leq 0.
\end{array}\right.
$$

Letting $n\to\infty$, we find
$$
\int_0^L \psi(x)(u(x)T(u(x))-\overline{u}(x)T(\overline{u}(x)))
{\rm sign}_0^+(u(x)-\overline{u}(x))\, dx
$$
$$
+\int_0^L{\rm
sign}_0^+(u(x)-\overline{u}(x))(T(u(x))\z(x)-T(\overline u(x))\zbar(x))\cdot \psi'\left(x\right) \, dx
$$
$$
\leq \int_0^L\psi(x)
[f(x)T(u(x))-\overline{f}(x)T(\overline{u}(x))] {\rm
sign}_0^+(u(x)-\overline{u}(x))\, dx + \int_0^L
\psi(x) \overline{\z} \, DT(\overline{u}).
$$

Taking now a sequence $\psi_m\uparrow \1_{]0,L[}$, $\psi_m\in
\mathcal D(]0,L[)$ in the above formula, we have
$$
\int_0^L (u(x)T(u(x))-\overline{u}(x)T(\overline{u}(x))) {\rm
sign}_0^+(u(x)-\overline{u}(x))\,
dx$$$$+\lim_{m\to\infty}\int_0^L{\rm
sign}_0^+(u(x)-\overline{u}(x))(T(u(x))\z(x)-T(\overline u(x))\zbar(x))\cdot \psi_m'\left(x\right) \, dx$$
$$\leq \int_0^L
(f(x)T(u(x))-\overline{f}(x)T(\overline{u}(x))) {\rm
sign}_0^+(u(x)-\overline{u}(x))\, dx + \int_0^L \overline{\z} \, DT(\overline{u}).
$$

Now we deal with the second term in the above expression.
$$
\lim_{m\to\infty}\int_0^L{\rm
sign}_0^+(u(x)-\overline{u}(x))(T(u(x))\z(x)-T(\overline u(x))\zbar(x))\cdot \psi_m'\left(x\right) \, dx
$$
$$
=-\lim_{m\to\infty}\int_0^L\psi_m(x) \{
  \z D[{\rm
sign}_0^+(u-\overline{u})T(u)] -  \overline{\z}(x) D[{\rm
sign}_0^+(u-\overline{u})T(\overline{u})]\}
$$
$$
+\lim_{m\to\infty}\int_0^L\psi_m(x) \{ {\rm
sign}_0^+(u(x)-\overline{u}(x))T(\overline{u}(x)) \overline{\z}'(x)
-{\rm
sign}_0^+(u(x)-\overline{u}(x))T(u(x)) \z'(x)\,
 \}\,
dx \, ,
$$
which leads to
$$
=-\int_0^L{\rm
sign}_0^+(u(x)-\overline{u}(x))T(u(x))\z'(x) \, dx -\int_0^L \z D[{\rm
sign}_0^+(u-\overline{u})T(u)]
$$
$$
+\int_0^L{\rm
sign}_0^+(u(x)-\overline{u}(x))T(\overline{u}(x)) \overline{\z}'(x) \,
dx +\int_0^L \overline{\z} D[{\rm
sign}_0^+(u-\overline{u})T(\overline{u})]
$$
$$
= \left[\z(0) T(u(0_{+}))-\overline{\z}(0)T(\overline{u}(0_{+})) \right] {\rm
sign}_0^+(u(0_{+})-\overline{u}(0_{+}))
$$
 $$
 -\left[\z(L)
T(u(L_{-}))-\overline{\z}(L)T(\overline{u}(L_{-})) \right] {\rm
sign}_0^+(u(L_{-})-\overline{u}(L_{-})).
$$
Therefore,
$$
\int_0^L (u(x)T(u(x))-\overline{u}(x)T(\overline{u}(x))) {\rm
sign}_0^+(u(x)-\overline{u}(x))\, dx$$
$$+ \left[\z(0) T(u(0_{+}))-\overline{\z}(0)T(\overline{u}(0_{+})) \right] {\rm
sign}_0^+(u(0_{+})-\overline{u}(0_{+})) $$ $$-\left[\z(L)
T(u(L_{-}))-\overline{\z}(L)T(\overline{u}(L_{-})) \right] {\rm
sign}_0^+(u(L_{-})-\overline{u}(L_{-}))$$ $$\leq \int_0^L
[f(x)T(u(x))-\overline{f}(x)T(\overline{u}(x))] {\rm
sign}_0^+(u(x)-\overline{u}(x))\, dx + \int_0^L \overline{\z} \, DT(\overline{u}).
$$

\noindent Dividing by $b > 0$, and letting $a\to 0^+$,
 and letting $b\to 0^+$ in this order, we obtain
$$
\begin{array}{l}
\displaystyle \int_0^L
(u\1_{[u>0]}-\overline{u}\1_{[\overline{u}>0]}) {\rm
sign}_0^+(u-\overline{u})\, dx \\ \\+ \left[\z(0) {\rm
sign}_0^+(u(0_{+}))-\overline{\z}(0){\rm
sign}_0^+(\overline{u}(0_{+})) \right] {\rm
sign}_0^+(u(0_{+})-\overline{u}(0_{+})) \\ \\ - \left[\z(L) {\rm
sign}_0^+(u(L_{-}))-\overline{\z}(L){\rm
sign}_0^+(\overline{u}(L_{-})) \right] {\rm
sign}_0^+(u(L_{-})-\overline{u}(L_{-}))
\\ \\
\leq \displaystyle\int_0^L (f
\1_{[u>0]}-\overline{f}\1_{[\overline{u}>0]}) {\rm
sign}_0^+(u-\overline{u})\, dx + \lim_{b\to
0}\frac{1}{b}\left(\lim_{a\to 0}\int_0^L \overline{\z} \,
DT(\overline{u})\right).
\end{array}
$$
Now, since $\z(0) = \overline{\z}(0) = - \beta \not=0$, and
$u(0_{+}) \geq \frac{\beta}{c} > 0$ and $\overline{u}(0_{+}) \geq
\frac{\beta}{c} > 0$  by (\ref{uencero}), we have  that the second
term in the above expression vanishes. On the other hand, since
$\z(L)= - c u(L_{-})$ and $\overline{\z}(L) =- c
\overline{u}(L_{-})$, the third  term in the above expression is
non-negative. Consequently,
\begin{eqnarray} \label{E1Des}
\displaystyle \int_0^L
(u\1_{[u>0]}-\overline{u}\1_{[\overline{u}>0]}) {\rm
sign}_0^+(u-\overline{u})\, dx \hspace{2.1cm} \nonumber
\\
\leq \displaystyle\int_0^L (f
\1_{[u>0]}-\overline{f}\1_{[\overline{u}>0]}) {\rm
sign}_0^+(u-\overline{u})\, dx +  \displaystyle\lim_{b\to
0}\frac{1}{b}\left(\lim_{a\to 0}\int_0^L \overline{\z} \,
DT(\overline{u})\right). \hspace{0.4cm}
\end{eqnarray}

 Next we claim
that
\begin{equation}\label{u=0} f = 0 \ \ {\rm a.e.
\ on} \ [u=0] \ \ \ \ {\rm and} \ \ \ \ \overline{f} = 0 \ \ {\rm
a.e \ on} \ [\overline{u} = 0].
\end{equation}
Let $0 \leq \phi \in {\mathcal D}(]0,L[)$ be and $a > 0, \epsilon
> 0 $. Multiplying $f - u = -\z'$ in ${\mathcal D}^{\prime}(]0,L[)$ by
$T_{a, a + \epsilon}^a(u)\phi$ and integrating by parts and having
in mind (\ref{Cond2,2})  and (\ref{positivity}), we have
$$
\int_0^L (f - u) T_{a, a + \epsilon}^a(u)\phi \, dx = \int_0^L\phi
\z DT_{a, a + \epsilon}^a(u) + \int_0^L \z \cdot  \phi' T_{a, a +
\epsilon}^a(u) \, dx \geq \int_0^L \z \cdot  \phi' T_{a, a +
\epsilon}^a(u) \, dx.
$$
Dividing  by $ \epsilon$ and letting $\epsilon\to 0^+$, we get
$$\int_0^L (f- u)\1_{[u > a]} \phi \, dx \geq \int_0^L
\z \cdot  \phi' \1_{[u > a]} \, dx.$$ Hence
$$\int_0^L (f - u)\1_{[u \leq a]} \phi \, dx = \int_0^L (f - u) \phi \,
dx- \int_0^L (f - u)\1_{[u > a]}(x) \phi \, dx$$
$$\leq \int_0^L (f - u) \phi \, dx - \int_0^L
\z \cdot \phi' \1_{[u > a]} \, dx = \int_0^L \z \cdot \phi' \1_{[u
\leq a]} \, dx.$$ Then, letting $a \to 0^+$, since $\z = 0$ in $[u
= 0]$, we have
$$ \int_0^L f\1_{[u = 0]} \phi \, dx = \int_0^L (f - u)\1_{[u = 0]} \phi \,
dx \leq 0,$$ for all $0 \leq \phi \in {\mathcal D}(]0,L[)$, from
where it follows that $f \1_{[u = 0 ]} = 0$ a.e. in $]0,L[$.
Similarly, $\overline{f} \1_{[\overline{u} = 0 ]} = 0$ a.e. in
$]0,L[$ and (\ref{u=0}) holds.

On the other hand, by (\ref{u=0})  we have,
$$\lim_{b\to
0}\frac{1}{b}\left(\lim_{a\to 0}\int_0^L \zbar DT(\overline
u)\right) $$ $$=-\lim_{b\to 0}\frac{1}{b}\lim_{a\to
0}\left(\overline\z(0)T(\overline{u}(0_{+}))- \overline\z(L)
T(\overline{u}(L_{-}))+\int_0^L T(\overline u)
\overline\z'\,dx\right)$$$$=-\lim_{b\to 0}\frac{1}{b}\left(
\overline\z(0)T_{0,b}(\overline{u}(0_{+})) -  \overline\z(L)
T_{0,b}(\overline{u}(L_{-}))+ \int_0^L T_{0,b}(\overline
u)\overline\z' \,dx \right)$$
$$=- \overline\z(0){\rm
sign}_0^+(\overline{u}(0_{+})) + \overline\z(L){\rm
sign}_0^+(\overline{u}(L_{-}))- \int_0^L \1_{[\overline
u>0]}\overline\z' \,dx $$ $$=- \overline\z(0){\rm
sign}_0^+(\overline{u}(0_{+}))  + \overline\z(L){\rm
sign}_0^+(\overline{u}(L_{-})) - \int_0^L \overline\z' \,dx  $$
$$= \overline\z(0) \left( 1 - {\rm
sign}_0^+(\overline{u}(0_{+}))  \right) + \overline\z(L)
\left({\rm sign}_0^+(\overline{u}(L_{-})) - 1 \right) = 0.$$

\noindent
Then, from (\ref{E1Des}), it
follows that
$$\displaystyle \int_0^L
(u\1_{[u>0]}-\overline{u}\1_{[\overline{u}>0]}) {\rm
sign}_0^+(u-\overline{u})\, dx \leq \displaystyle\int_0^L
(f\1_{[u>0]}-\overline{f}\1_{[\overline{u}>0]}) {\rm
sign}_0^+(u-\overline{u})\, dx.$$ Hence, using (\ref{u=0}), we
obtain
$$
\int_0^L ( u-\overline{u})^+ \, dx \leq \int_0^L (f-\overline{f})
{\rm sign}_0^+(u-\overline{u})\, dx \leq \int_0^L
(f-\overline{f})^+ \, dx.
$$
This concludes the proof of the uniqueness part of Theorem
\ref{EUTElliptic}.\hfill$\Box$

\section{Semigroup solution}\label{SS}

In this section we shall associate an accretive operator  in $L^1(]0, L[)$ to the problem (\ref{Elliptictproblem}).

 \begin{definition}{\rm $(u, v) \in  \mathcal{B}_{\beta}$ if and only
 if $0 \leq u \in TBV^+(]0, L[)$,
$ v \in L^{1}(]0, L[)$ and $u$ is the entropy solution of problem
(\ref{Elliptictproblem})}.
\end{definition}

From Theorem \ref{EUTElliptic}, it follows that the operator
$\mathcal{B}_{\beta}$ is $T$-accretive in $L^1(]0, L[)$ and
verifies
\begin{equation}\label{Rcond00}
 L^{\infty}(]0, L[)^+
\subset R(I + \lambda {\mathcal B}_{\beta}) \quad \hbox{for all} \
\ \lambda
> 0.
\end{equation}
In order to get an $L^\infty$-estimate of the resolvent, we need to find the  steady
state solution, that is, the
 function $u_{\beta}$ which is the entropy solution of the problem
 \begin{equation}\label{ElliptictproblemSS}\left\{
\begin{array}{ll}
  - \left(\a(u_{\beta},u_{\beta}') \right)'  =  0 & \ \hbox{in} \ \  ]0, L[
\\ \\ -  \a(u_{\beta},u_{\beta}')|_{x=0}  = \beta > 0 \ \ \hbox{and} \ \ u_{\beta}(L)=0.
\end{array}
\right.
\end{equation}

\begin{proposition}\label{Cojonudo2}  There is a non-increasing function
$u_{\beta} \in C^1(]0, L[)$,  with $u_{\beta} \geq
\frac{\beta}{c}$, that is an entropy solution of the stationary
problem (\ref{ElliptictproblemSS}). Moreover, there exists a
constant  $M:= M(c, \beta, \nu, L)$ such that
$$ \Vert u_{\beta} \Vert_{\infty} \le M.$$
\end{proposition}\noindent{\bf Proof.}

Integrating (\ref{ElliptictproblemSS}) over $]0,L[$ we find that
$\a(u_{\beta},u_{\beta}')(L)  =  - \beta$. Now, if $u_\beta$ has to fulfil the weak
Dirichlet condition $\a(u_{\beta},u_{\beta}')(L)  = -c u_\beta (L_-)$  then we must have
$u_\beta(L_-) = \beta /c$. We will follow this prescription hereafter.

\noindent
If $u_{\beta}$ is a solution of the problem
(\ref{ElliptictproblemSS}), we have
$$- \left(\a(u_{\beta},u_{\beta}') \right)'  =  0 \iff \nu
\frac{u_{\beta} u_{\beta}'}{\sqrt{u_{\beta}^2 + \frac{\nu^2}{c^2} (u_{\beta}')^2}} = -
\beta.$$ Then, assuming that $u_{\beta}' < 0$, we get
$$u_{\beta}' = - \frac{\beta \, u_{\beta}}{ \nu \sqrt{{u_{\beta}^2
- \left( \frac{\beta}{c} \right)^2}}}.$$ Thus, we get that
$u_{\beta}$ satisfies the ordinary differential equation
$$
\frac{u_{\beta}' \sqrt{{u_{\beta}^2 - \left(\frac{\beta}{c}
\right)^2}}}{u_{\beta}} = - \frac{\beta}{\nu}.
$$
By means of the change of variable $v^2 = u_{\beta}^2 -
\left(\frac{\beta}{c} \right)^2$, we arrive to the ODE
$$
- \frac{\beta}{\nu} = \left( 1 - \frac{1}{1 + \left( \frac{v}{ \beta /c} \right)^2}
\right) v'.
$$
Then,
$$
\int_x^L \left(- \frac{\beta}{\nu}\right) \, dy = \int_x^L v'(y)
\, dy - \int_x^L \frac{v'(y)}{1 + \left(\frac{v(y)}{ \beta
/c}\right)^2} \, dy $$ $$= v(L) - v(x) - \frac{\beta}{c} \arctan
\left( \frac{v(L)}{ \beta /c} \right) + \frac{\beta}{c} \arctan
\left( \frac{v(x)}{\beta /c}\right).
$$
Hence, we get
\begin{equation}\label{ODE2}
x = L - \frac{\nu}{\beta} \sqrt{u_{\beta}(x)^2 - \left(\frac{\beta}{c} \right)^2}+
\frac{\nu}{c} \arctan \left[\frac{c}{\beta}\sqrt{u_{\beta}(x)^2 - \left(\frac{\beta}{c}
\right)^2} \right].
\end{equation}
 If $x = u_{\beta}^{-1}(y)$, then
we can write (\ref{ODE2}) as
$$u_{\beta}^{-1}(y)  = L -
\frac{\nu}{\beta}
 \sqrt{y^2 - \left(\frac{\beta}{c}
\right)^2}+ \frac{\nu}{c} \arctan \left[\frac{c}{\beta}\sqrt{y^2 -
\left(\frac{\beta}{c} \right)^2} \right].
$$
Thus,
$$\left(u_{\beta}^{-1}\right)^{\prime}(y) = \frac{y}{\sqrt{y^2 - \left(\frac{\beta}{c}
\right)^2}} \left(\frac{\nu}{\beta} \right) \left( -1 +
\frac{\beta^2}{c^2 y^2}\right),
$$
and consequently, since $(u_{\beta})(L_{-}) = \frac{\beta}{c}$, we
obtain that
$$
(u_{\beta})^{\prime}(L_{-}) = \lim_{y \searrow \frac{\beta}{c}}
\frac{1}{\left(u_{\beta}^{-1}\right)^{\prime}(y)} = \lim_{y \searrow
\frac{\beta}{c}} \frac{\sqrt{y^2 - \left(\frac{\beta}{c}
\right)^2}}{y} \left(\frac{\beta}{\nu} \right) \left( \frac{c^2
y^2}{\beta^2 - c^2 y^2}\right) = - \infty.
$$
Finally, since $u_\beta$ satisfies $-(\a(u_{\beta}(x), u_{\beta}'(x)))' = 0$ if $x \in ]0,L[$ and satisfies the boundary conditions also, we have that $u_\beta$ is an entropy solution of the problem (\ref{ElliptictproblemSS}).

\hfill$\Box$

The following homogeneity of the operator ${\mathcal B}_{\beta}$ will be important to get
the $L^{\infty}$-estimate of the resolvent.

\begin{proposition}\label{Cojonudo1} For $\mu > 0, \lambda > 0$ and $\beta > 0$, we
have
\begin{equation}\label{homog1}\left(I + \lambda {\mathcal B }_{\beta} \right)^{-1} (\mu u) =
\mu \left(I + \lambda {\mathcal B }_{\frac{\beta}{\mu}}
\right)^{-1} (u).
\end{equation}
Moreover, for $\beta_ 1 \leq \beta_2$,  $u \in  L^{\infty}(]0, L[)^+$ and $\lambda > 0$ such that $\left(I + \lambda {\mathcal B }_{\beta_2} \right)^{-1} (u) \in BV(]0,L[)$, we have
\begin{equation}\label{homog3}
\left(I +  \lambda {\mathcal B }_{\beta_1} \right)^{-1} (u) \leq
\left(I + \lambda {\mathcal B }_{\beta_2} \right)^{-1} (u) \ \ \hbox{a.e.}\ x\in ]0,L[.
\end{equation}

\end{proposition}
\noindent{\bf Proof.} From the definition of the operator it is
easy to see that if $u \in D({\mathcal B }_{\frac{\beta}{\mu}})$,
then $\mu u \in D({\mathcal B}_{\beta})$ and ${\mathcal B}_{\beta}(\mu u) = \mu {\mathcal B
}_{\frac{\beta}{\mu}}(u).
$
Then, we have
$$v := \left(I + \lambda {\mathcal B }_{\beta} \right)^{-1} (\mu
u) \iff v + \lambda {\mathcal B}_{\beta} (v) = \mu u \iff
\frac{1}{\mu} v + \frac{1}{\mu} \lambda {\mathcal B }_{\mu
\frac{\beta}{\mu}}(v) = u $$ $$\iff  \frac{1}{\mu} v + \lambda
{\mathcal B }_{\frac{\beta}{\mu}}( \frac{v}{\mu}) = u \iff \left(I
+ \lambda {\mathcal B }_{\frac{\beta}{\mu}} \right)^{-1} (u) =
\frac{v}{\mu},$$ from where (\ref{homog1}) follows.

Finally, let us see that (\ref{homog3}) holds. Let $u_i:= \left(I
+  \lambda{\mathcal B }_{\beta_i} \right)^{-1} (u)$, $i = 1,2.$
Then, $u_i$ is an entropy solution of the problem

\begin{equation}\nonumber\label{ElliptictproblemYY}\left\{
\begin{array}{ll}
 u_i - \lambda \left(\a(u_i,u_i') \right)'  =  u & \ \hbox{in} \ \  ]0, L[
\\ \\ -  \a(u_i,u_i')|_{x=0}  = \beta_i > 0 \ \ \hbox{and} \ \ u(L)=0.
\end{array}
\right.
\end{equation}
Therefore, if $p_n$ are non  negative increasing functions that are an approximation of
the ${\rm sign}_0^+$ function, having in mind (\ref{monotoneineq}), since $p_n(u_1 - u_2) \in BV(]0,L[)$,  we  get
$$\int_0^L (u_1- u_2) p_n(u_1 - u_2) \, dx = \int_0^L \lambda \left(\left(\a(u_1,u_1')
\right)'- \left(\a(u_2,u_2') \right)' \right)  p_n(u_1 - u_2) \,
dx
$$ $$= - \int_0^L \lambda \left(\a(u_1,u_1') - \a(u_2,u_2') \right) D( p_n(u_1 - u_2))
$$ $$+
 \lambda \left( \a(u_1,u_1')(L_-) - \a(u_2,u_2')(L_-) \right) p_n(u_1 - u_2) (L_{-})
-  \lambda \left( \a(u_1,u_1')(0_+) - \a(u_2,u_2')(0_+) \right) p_n(u_1
- u_2) (0_{+})$$
 $$\leq
 \lambda \left( \a(u_1,u_1')(L_-) - \a(u_2,u_2')(L_-)\right) p_n(u_1 - u_2) (L_{-})
+ \lambda(\beta_1 - \beta_2) p_n(u_1 - u_2) (0_{+})$$ $$\leq
 \lambda \left( \a(u_1,u_1')(L_-) - \a(u_2,u_2')(L_-) \right) p_n(u_1 - u_2)
 (L_{-}).$$
 Then, taking limit as $n \to  + \infty$ we get
 $$\int_0^L (u_1- u_2)^+  \, dx \leq
 \lambda \left( \a(u_1,u_1')(L_-) - \a(u_2,u_2')(L_-) \right) \hbox{sign}^+_0(u_1 - u_2)
 (L_{-}) \leq 0,$$
 since $\a(u_i, u_i')(L) = - c u_i(L_{-}) $, $i = 1,2.$
Therefore, $u_1 \leq u_2$, and we
 finish the proof.
 \hfill$\Box$

\begin{proposition}\label{Cojonudo99} For $u \in  L^{\infty}(]0, L[)^+$ and $\lambda >
0$, we have
\begin{equation}\nonumber\label{homog39}
0 \leq \left(I +  \lambda {\mathcal B }_{\beta} \right)^{-1} (u)
\leq \mu u_{\beta}, \ \ \ {\rm with} \ \ \mu= \max \left\{\frac{c
\Vert u \Vert_{\infty}}{\beta}, 1 \right\}.
\end{equation}

\end{proposition}
\noindent{\bf Proof.} Let $u_{\beta}$ be the  entropy solution of
the stationary problem (\ref{ElliptictproblemSS}) given in
Proposition \ref{Cojonudo2}. Then, $(u_{\beta}, 0) \in {\mathcal B
}_{\beta}$, from where it follows that
\begin{equation}\label{0Barsa1}
\left(I + \lambda{\mathcal B }_{\beta} \right)^{-1}(u_{\beta}) =
u_{\beta}.
\end{equation}

On the other hand, since $u_{\beta} \geq \frac{\beta}{c}$, if $\mu:= \max \{ \frac{c
\Vert  u \Vert_{\infty}}{\beta}, 1 \}$, we have $0 \leq  u \leq \mu u_{\beta}$. Hence, by
Proposition \ref{Cojonudo1} and having in mind (\ref{0Barsa1}), we get
$$0 \leq  \left(I + \lambda  {\mathcal B }_{\beta} \right)^{-1}(u)
\leq \left(I + \lambda {\mathcal B }_{\beta} \right)^{-1}(\mu
u_{\beta}) = \mu \left(I + \lambda {\mathcal B
}_{\frac{\beta}{\mu}} \right)^{-1} (u_{\beta})\leq  \mu \left(I +
\lambda {\mathcal B }_{\beta} \right)^{-1}(u_{\beta}) = \mu
 u_{\beta}.$$

\hfill$\Box$

Next we introduce the main result of this section, which paves the way for the operator $\mathcal{B}_{\beta}$ to generate an order-preserving semigroup \cite{BCrP}.

\begin{theorem}
\label{condrango}
 $\mathcal{B}_{\beta}$ is $T$-accretive in $L^1(]0, L[)$, and verifies the range condition
\begin{equation}\nonumber\label{Rcond}
\overline{D({\mathcal B}_{\beta})}^{L^1(]0, L[)} = L^1(]0, L[)^+
\subset R(I + \lambda {\mathcal B}_{\beta}) \quad \hbox{for all} \
\ \lambda
> 0.
\end{equation}

\end{theorem}

\noindent{\bf Proof.} The $T$-accretivity of the operator
$\mathcal{B}_{\beta}$ is known, and that it verifies (\ref{Rcond00}) also. To
prove the density of $D(\mathcal{B}_{\beta})$ in $L^1(]0, L[)^+$,
we prove that ${\mathcal D}(]0, L[)^+ \subseteq
\overline{D(\mathcal{B}_{\beta})}^{L^1(]0, L[)}$. Let $0 \leq v
\in {\mathcal D}(]0, L[)$. By (\ref{Rcond00}), $v \in R(I +
\frac{1}{n}
 \mathcal{B}_{\beta})$ for all $n \in \N$. Thus, for each $n \in \N$, there
exists $u_{n} \in D(\mathcal{B})$ such that $(u_{n}, n(v-u_{n}))
\in \mathcal{B}$. Since $u_n = (I + \frac{1}{n}
 \mathcal{B}_{\beta})^{-1}(v)$, by Proposition \ref{Cojonudo99},
 we get
 \begin{equation}\label{ACOTTA}\Vert u_n\Vert_\infty \leq M:= M(\beta, c, \nu, L, \Vert v
\Vert_\infty).
\end{equation}
Let $\epsilon >0$. Since
$$n(v - u_n) = - D\a(u_n,
 u_n^{\prime}) \ \ \ \ {\rm in} \ \ \ \ {\mathcal
D}^{\prime}(]0,L[),$$ multiplying by $v - S_{\epsilon}(u_n)$,
with $S_{\epsilon}:= T_{\epsilon, \Vert v \Vert_{\infty}}$, and
integrating by parts, we get
$$\displaystyle \int_0^L (v - S_{\epsilon}(u_n)) n(v-u_{n}) \ dx = \left[ \int_0^L \a(u_n,
 u_n^{\prime})(Dv  - DS_{\epsilon}(u_n))\right]$$ $$ - c u_n(L_{-})S_{\epsilon}(u_n)(L_{-})
  + \beta S_{\epsilon}(u_n)(0_{+}).
  $$
Then, since
 $$\int_0^L \a (u_n,u_n') D S_\epsilon (u_n) \ge 0,$$
having in mind (\ref{ACOTTA}), we get
$$\displaystyle \int_0^L (v - S_{\epsilon}(u_n)) (v-u_{n}) \ dx \leq \frac{1}{n} \left[ \int_0^L \a(u_n,
 u_n^{\prime})Dv\right]
 + \frac{1}{n} \beta S_{\epsilon}(u_n)(0_{+}) \leq \frac{C}{n}.$$
Letting $\epsilon \to 0^+$, we get $$\int_0^L(v-u_{n})^2 \ dx \leq
\frac{C}{n},$$ and we obtain that $u_{n} \rightarrow v$ in
$L^2(]0, L[)$, as $n \rightarrow \infty$. Moreover, we have $u_{n}
\rightarrow v$ in $L^1(]0, L[)$, as $n \rightarrow \infty$.
Therefore $v \in \overline{D(\mathcal{B}_{\beta})}^{L^1(]0, L[)}$
and the proof of the density of $D(\mathcal{B}_{\beta})$ in
$L^1(]0, L[)^+$ is complete.

To finish the proof of the theorem, we only need to show that the
operator $\mathcal{B}_{\beta}$ is closed in $L^1(]0,L[) \times
L^1(]0,L[)$. Given $(u_n, v_n) \in \mathcal{B}_{\beta}$ such that
$u_n \to u$ and $v_n \to v$ in $L^1(]0,L[)$, we need to prove that
$(u, v) \in \mathcal{B}_{\beta}$. Since $(u_n, v_n) \in
\mathcal{B}_{\beta}$, we have that $u_n \in TBV^+(]0,L[)$ and
$\z_n:= \a(u_n,
 u_n') \in C([0,L])$ satisfy
\begin{equation}\label{TTCond1}
 v_n  = - D \z_n\ \ \ \ \ \ {\rm in} \ \ {\mathcal
D}^{\prime}(]0,L[),
\end{equation}
\begin{equation}\label{TTCond2,2} h(u_n,DT(u_n))\le \z_nDT(u_n)
\ \ \ {\rm as \ measures} \ \ \forall \, T \in {\mathcal T}^+
\end{equation}
\begin{equation}\nonumber\label{TTSCond2}
\displaystyle   h_S(u_n,DT(u_n))\le \z_nDJ_{T'S}(u_n) \ \ \ \hbox{as  measures}
\\ \\  \ \displaystyle\forall \, S \in {\mathcal P}^+,  \, T\in {\mathcal T}^+,
\end{equation}
\begin{equation}
\label{TTCond2}
- \z_n(0) =  \beta \ \ \  {\rm and} \ \ \ \
  \z_n(L) = - c u_n(L_{-}).
\end{equation}

\noindent
Let $T = T_{a,b} \in \mathcal{T}_r$. Multiplying (\ref{TTCond1}) by $T(u_n)$ and
applying integration by parts (Lemma \ref{IntBP}), we get
$$\int_0^L v_n T(u_n) \, dx  = \int_0^L \z_n DT(u_n) - \z_n(L)T(u_n(L_-)) - \beta T(u_n(0_+)),$$
from where it follows that
\begin{equation}\label{TTEE1}
\int_0^L \z_n DT(u_n) \leq b(\beta + \Vert v \Vert_1) \le C.
\end{equation}
Here we used the boundary condition \eqref{TTCond2} to be able to disregard the term related to
 $\z_n(L)$, as it has the right sign.

On the other hand, by (\ref{TTCond2,2}) and having in mind (\ref{BOUNFG}), we get
\begin{equation}\label{TTEE2}
\int_0^L \z_n DT(u_n) \geq  \frac{c}{2} \int_0^L \vert D([T(u_n)]^2) \vert  - \frac{c^2}{\nu}
\int_0^L T(u_n)^2 \, dx.
\end{equation}
By (\ref{TTEE1}) and (\ref{TTEE2}), we obtain that
\begin{equation}
\label{TTEE3}
  \int_0^L \vert D([T(u_n)]^2) \vert \leq \frac{2c}{\nu}
\int_0^L T(u_n)^2 \, dx + \frac{2  C}{c} \leq  \frac{2cLb^2}{\nu} + \frac{2 C}{c} =
 C.
\end{equation}
Using the coarea formula as in the proof of Theorem \ref{EUTElliptic}, from (\ref{TTEE3})
we deduce that
$$
 \int_0^L \vert DT(u_n) \vert \leq \frac{C}{2a}
 \quad \quad \forall \, n \in \N.
$$
Then, since the total variation is semi-continuous in $L^1(]0,L[)$, we have
\begin{equation}\nonumber\label{TTEE4}
\int_0^L \vert DT(u) \vert  \leq \liminf_{n \to \infty} \int_0^L \vert DT(u_n) \vert \leq
\frac{C}{2a}.
\end{equation}
Hence, $T(u) \in BV(]0,L[)$, and consequently, $u \in TBV^+(]0,L[)$.

Since $\z_n = c \vert u_n \vert \b(u_n, u^{\prime}_n)$ with $\vert \b(u_n, u^{\prime}_n)
\vert \leq 1$, for all measurable subsets $E \subset ]0,L[$, we have
$$
\int_E \vert \z_n \vert \, dx \leq c \int_E \vert u_n \vert \, dx.
$$
Therefore, by Dunford-Pettis's Theorem, we can assume that
\begin{equation}\label{TTEE5}
\z_n \rightharpoonup \z \ \ \ \hbox{weakly in} \ L^1(]0,L[).
\end{equation}
Moreover, since $\vert \b(u_n, u^{\prime}_n) \vert \leq 1$, we also can assume that
\begin{equation}\label{TTEE6}
\b(u_n, u^{\prime}_n) \rightharpoonup \z_b \ \ \ \hbox{weakly}^* \ \hbox{in} \
L^{\infty}(]0,L[).
\end{equation}
As $u_n \to u$ in $L^1(]0,L[)$, from (\ref{TTEE5}) and (\ref{TTEE6}), we obtain that
\begin{equation}\label{TTEE7}
\z = c u \z_b.
\end{equation}
As $v_n \to v$ in $L^1(]0,L[)$, from (\ref{TTEE5}) and (\ref{TTCond1}), we easily deduce
that
\begin{equation}\label{Cond1098}
 v  = - D \z \ \ \ \ \ \ {\rm in} \ \ {\mathcal
D}^{\prime}(]0,L[),
\end{equation}
and by (\ref{TTEE7}) and (\ref{Cond1098}), we have $\z \in W^{1,1}(]0,L[) \subset
C([0,L])$.

\begin{lemma}\label{MintBrrN}
\begin{equation}\nonumber\label{ident_campo1}
 \z(x)=\a(u(x),u^\prime(x)) \ \ \mbox{a.e.} \, \, x \in ]0,L[
\end{equation}
\end{lemma}
\noindent{\bf Proof.}
We use Minty-Browder's technique. Let
$0<a<b$, let $0 \le \phi \in C_c^1(]0,L[)$ and let $g \in C^2([0,L])$. By
(\ref{monotoneineq}), we have that
\begin{equation}
\label{monotono}
   \int_0^L \phi [\a(u_n,u_n^\prime)-\a(u_n,g^\prime)] T_{a,b}^\prime(u_n)(u_n-g)^\prime \, dx \ge 0.
\end{equation}
Let us denote
$$
    J_{\a}(x,r) := \int_0^r \a(s,g^\prime(x)) \, ds,
$$
$$
   J_{\a^\prime}(x,r):= \int_0^r \partial_x[\a(s,g^\prime(x))] \, ds = \int_0^r \frac{\partial \a}{\partial \xi}(s,g^\prime(x)) g^{\prime \prime}(x) \, ds
$$
and observe that
$$
   - \a\left(T_{a,b}(u_n(x)),g^\prime(x) \right) [T_{a,b}(u_n)]'
=  -D^{ac} \left[J_\a(x,T_{a,b}(u_n(x)))\right]  +
J_{\a^\prime}(x,T_{a,b}(u_n(x))),
$$
this we will substitute into (\ref{monotono}).
Note now that, using (\ref{TTCond2,2})
$$
  \int_0^L \phi \left[\z_n D^sT_{a,b}(u_n) - D^s J_{\a} (x,T_{a,b}(u_n)) \right]
  \ge \int_0^L \phi \left[h(u_n,DT_{a,b}(u_n))^s - D^s J_{\a} (x,T_{a,b}(u_n)) \right] \ge 0,
$$
where the last inequality is proved using the properties of the Lagrangian (see \cite{ARMA}). Then we can add this inequality to (\ref{monotono}):
$$
0 \le \int_0^L \phi \left[ \z_n DT(u_n) - DJ_\a (x,T_{a,b}(u_n(x))) \right]
$$
$$
+ \int_0^L \phi \left[J_{\a'}(x,T(u_n(x)))  - \z_n g' T_{a,b}^\prime(u_n)) + g' T_{a,b}^\prime(u_n) \a(u_n,g') \right]\ dx.
$$
Now, since
$$
    \int_0^L \phi \, \z_n [D T_{a,b}(u_n) - g' T_{a,b}^\prime(u_n)]
  =    \int_0^L \phi \, \z_n D [T_{a,b}(u_n) - g ] +    \int_0^L \phi \, \z_n g' (1-T_{a,b}^\prime(u_n)) \, dx
$$
$$
 = - \int_0^L v_n \phi \left(T_{a,b}(u_n) - g \right) \, dx
   - \int_0^L \left(T_{a,b}(u_n) - g  \right) \a(u_n,u_n^\prime) \phi^\prime \, dx +  \int_0^L \phi \, \z_n g' (1-T_{a,b}^\prime(u_n)) \, dx
$$
we get
$$
  \lim_{n \to + \infty} \int_0^L \phi \, \z_n [D T_{a,b}(u_n) - g' T_{a,b}^\prime(u_n)]\, dx
  \le \langle \z D (T_{a,b}(u)-g), \phi \rangle + \|g^\prime\|_\infty \int_0^L \vert \z \vert \phi \left(1- T_{a,b}^\prime(u) \right) \, dx.
$$
On the other hand, the almost everywhere convergence of $u_n$ implies that
$$
   J_{\a^\prime}\left(x,T_{a,b}(u_n(x)) \right) \rightarrow J_{\a^\prime}\left(x, T_{a,b}(u(x)) \right) \, \mbox{a.e.}
$$
and we also have (see \cite{Ambrosio}, Proposition 3.13) that
$$D \left[J_\a(x,T_{a,b}(u_n(x))) \right] \rightharpoonup D \left[J_\a(x,T_{a,b}(u(x)))
\right] \, \mbox{weakly as measures}.
$$
As a consequence, we have
$$
   \lim_{n \to + \infty} \int_0^L \phi \left[J_{\a'}(x,T_{a,b}(u_n(x))) - D J_{\a}(x,T_{a,b}(u_n(x))) + g' T_{a,b}^\prime(u_n) \a(u_n,g') \right]
$$
$$
  = \langle J_{\a'}(x,T_{a,b}(u)) - D J_{\a}(x,T(u), \phi \rangle + \int_0^L \phi g' \a(u,g') T_{a,b}'(u) \ dx.
$$

Consequently we obtain
$$
   \langle \z \,D\left(T_{a,b}(u) - g \right), \phi \rangle
  + \|g^\prime\|_\infty \int_0^L \vert \z \vert \phi \left(1- T_{a,b}'(u) \right) \, dx
$$
$$
  + \int_0^L \phi \, \a(u,g^\prime) g^\prime T_{a,b}^\prime(u) \, dx
- \langle D \left[J_{\a}(x,T_{a,b}(u(x)))\right] -
J_{\a^\prime}(x,T_{a,b}(u(x))), \phi
 \rangle
 \ge 0
$$
for all $0\le \phi \in C_c^1(]0,L[)$. This means that, as measures,
$$
   \z \,D\left(T_{a,b}(u) - g \right)  - D \left[J_a(x,T_{a,b}(u(x)))\right]
    + J_{\a^\prime}(x,T_{a,b}(u(x)))
$$
$$
   + \left\{\a(u,g^\prime) g^\prime T_{a,b}^\prime(u) +
    \vert \z \vert  \|g^\prime\|_\infty \left(1- T_{a,b}^\prime(u) \right) \right\} \mathcal{L}^1 \ge 0,
$$and we obtain
$$
    \z \left( T_{a,b}(u) - g \right)^\prime - \a(u,g^\prime)(T_{a,b}(u))^\prime
+ \a(u,g^\prime)\, g^\prime T_{a,b}^\prime(u) + \vert \z \vert  \|g^\prime\|_\infty
\left(1- T_{a,b}^\prime(u) \right) \ge 0.
$$
If $x \in [a<u<b]$, this reduces to
$$
   \big(\z - \a(u,g^\prime) \big) (u-g)^\prime \ge 0,
$$
which holds for all $x\in \Omega \cap [a<u<b]$, where $\mathcal{L}^1(]0,L[ \setminus
\Omega) = 0$, and all $g \in C^2([0,L])$. Being $x \in \Omega \cap [a<u<b] $ fixed and
$\xi \in \RR$ given, we find $g$ as above such that $g^\prime (x) = \xi$. Then
$$
    \left(\z(x) - \a(u(x),\xi)\right) (u^\prime(x) - \xi) \ge 0, \, \, \forall \xi \in \RR.
$$
By an application of Minty-Browder's method in $\RR$, these inequalities imply that
$$
   \z(x) = \a(u(x),u^\prime (x)) \quad \mbox{a.e. on}\, [a<u<b].
$$
Since this holds for any $0<a<b$, we obtain the identification a.e. on the points of
$]0,L[$ such that $u(x) \ne 0$. Now, by our assumptions on $\a$ and (\ref{TTEE7}) we
deduce that $\z(x) = \a(u(x),u^\prime(x)) = 0$ a.e. on $[u=0]$. The Lemma is proved.

\hfill$\Box$

To finish the proof we only need to show that
\begin{equation}\nonumber\label{DDCondmeasure1}
\frac{c}{2} \vert D^s(T(u)^2) \vert \leq \z D^sT(u) \ \ \ {\rm as \ measures} \ \ \forall
\, T \in {\mathcal T}^+,
\end{equation}
\begin{equation}\nonumber\label{DDSCond2}
\displaystyle \vert D^s(J_{S\theta}(T(u)) \vert \leq \z D^sJ_{T'S}(u)  \ \ \hbox{as
measures}, \ \ \forall \, S \in {\mathcal P}^+ , \, T\in {\mathcal T}^+,
\end{equation}
\begin{equation}\nonumber\label{DDCondmeasure}
 - \a(u, u')(0) =  \beta  \ \ {\rm
and} \ \ \ \a(u, u')(L)  = - c u(L_{-}).
\end{equation}
These proofs are similar to those in the previous section.

\hfill$\Box$

From  Theorem \ref{condrango},  according to Crandall-Liggett's Theorem (c.f., e.g.,
\cite{BCrP}), for any $0 \leq u_0 \in L^1(]0,L[)$ there exists a unique mild solution $u
\in C([0, T]; L^1(]0,L[))$ of the abstract Cauchy problem
\begin{equation}\nonumber\label{ACP}
u^{\prime} (t) + {\mathcal B}_{\beta}u(t) \ni 0, \ \ \ \ u(0) =
u_0.
\end{equation}
Moreover, $u(t) = T_{\beta}(t) u_0$ for all $t \geq 0$, where $(T_{\beta}(t))_{t \geq 0}$
is the semigroup in $L^1(]0,L[)^+$ generated by Crandall-Liggett's exponential formula,
i.e.,
$$T_{\beta}(t) u_0 = \lim_{n \to \infty} \left(I + \frac{t}{n} {\mathcal B}_{\beta}\right)^{-n} u_0.$$ On the
other hand, as the operator ${\mathcal B }_{\beta}$ is
$T$-accretive we have that the comparison principle also holds for
$T_{\beta}(t)$, i.e., if $u_0, \overline{u}_0 \in L^1(]0,L[)^+$,
we have the estimate
\begin{equation}
\label{CPsemigroups} \Vert (T_{\beta}(t)u_0 -
T_{\beta}(t)\overline{u}_0)^+ \Vert_1 \leq \Vert (u_0 -
\overline{u}_0)^+\Vert_1.
\end{equation}

Obviously,  by Crandall-Liggett's exponential formula, from
(\ref{homog1}), we get that for all $u_0 \in L^1(]0,L[)^+$,
\begin{equation}\label{homog2}
T_{\beta}(t)(\mu u_0) = \mu T_{\frac{\beta}{\mu}}(t)(u_0) \ \ \ \
\hbox{for all} \ t > 0.
\end{equation}

As a consequence of (\ref{CPsemigroups}) and (\ref{homog2}), for $u \in  L^{\infty}(]0,
L[)^+$, we have

\begin{equation}\nonumber\label{34homog39}
0 \leq T_{\beta}(t) (u) \leq \mu u_{\beta}, \ \ \ {\rm with} \ \ \mu= \max \left\{\frac{c
\Vert u \Vert_{\infty}}{\beta}, 1 \right\}, \ \ \ \forall \, t \geq 0.
\end{equation}

\section{Existence and uniqueness of solutions of the parabolic
problem}\label{ParabolicProblem} 

This section deals with the problem
\begin{equation}
\label{Dirichletproblem} \left\{
\begin{array}{ll}
\displaystyle \frac{\partial u}{\partial t} = \left(\a(u,u_x) \right)_x
 &
\hspace{0.2cm}\hbox{in \hspace{0.15cm} $]0,T[\times ]0,L[$}\\
\displaystyle
 \\
-  \a(u(t,0),u_x(t,0))  = \beta > 0  \hbox{ and }  u(t,L)=0&
\hspace{0.2cm} \hbox{on \hspace{0.15cm} $t \in  ]0, T[$,}
 \\
\displaystyle
 \\
\displaystyle u(0,x) = u_{0}(x) & \hspace{0.2cm} \hbox{in \hspace{0.15cm} $x \in ]0,L[$.}
\end{array}
\right.
\end{equation}
To make precise our notion of solution we need to recall the
following definitions given in
\cite{ACMBook}. We set $Q_T=]0,T[ \times
]0,L[$.

It is well known (see for instance \cite{Schwartz}) that the dual
space $\left[L^1(0, T; BV(]0,L[))\right]^*$ is isometric to the
space $L^{\infty}(0, T; BV(]0,L[)^*, BV(]0,L[))$ of all weakly$^*$
measurable functions $f : [0, T] \rightarrow BV(]0,L[)^*$, such
that  $v(f) \in L^{\infty}([0, T])$,  where $v(f)$ denotes the
supremum of the set $\{ \vert \langle w, f \rangle \vert \ : \
\Vert w \Vert_{BV(]0,L[)} \leq 1 \}$ in the vector lattice of
measurable real functions. Moreover, the duality pairing is
$$\langle w, f \rangle = \int_0^T \langle w(t), f(t) \rangle \ dt,$$
for $w \in L^1(0, T; BV(]0,L[))$ and $f \in L^{\infty}(0, T;
BV(]0,L[)^*, BV(]0,L[))$.

By $L^1_{w}(0,T,BV(]0,L[))$ we denote the space of weakly
measurable functions $w:[0,T] \to BV(]0,L[)$ (i.e., $t \in [0,T]
\to \langle w(t),\phi \rangle$ is measurable for every $\phi \in
BV(]0,L[)^*$) such that $\int_0^T \Vert w(t)\Vert \, dt< \infty$.
Observe that, since $BV(]0,L[)$ has a separable predual (see
\cite{Ambrosio}), it follows easily that the map $t \in [0,T]\to
\Vert w(t) \Vert$ is measurable. By  $L^1_{loc, w}(0, T,
BV(]0,L[))$ we denote the space of weakly measurable functions
$w:[0,T] \to BV(]0,L[)$ such that the map $t \in [0,T]\to \Vert
w(t) \Vert$ is in $L^1_{loc}(]0, T[)$.

Note that if $w \in L^1(0, T; BV(]0,L[)) \cap L^{\infty}(Q_T)$ and $\z \in L^1(Q_T)$ is
such that there exists an element $\xi \in [L^1(0, T; BV(]0,L[))]^*$ with $D_x \z =\xi$
in $\mathcal{D}'(Q_T)$, we define, associated with  $(\z, \xi)$, the distribution $\z
D_xw$ in $Q_T$ by
\begin{equation}
\label{DefMRAd}
   \langle  \z D_x w, \varphi \rangle = -  \langle \xi, \varphi w \rangle
   - \int_0^T \int_0^L  \z(t,x) w(t,x) \partial_x \varphi(t,x) \, dxdt
\end{equation}
for all $\varphi \in \mathcal{D}(Q_T)$.

Our concept of solution for the problem (\ref{Dirichletproblem}) is the following.

\begin{definition} \label{DefESevP} {\rm
 A measurable function $u: ]0,T[\times ]0,L[ \rightarrow  \R^+$ is an {\it entropy solution}
of (\ref{Dirichletproblem}) in $Q_T = ]0,T[\times ]0,L[$ if $u \in C([0, T];
L^1(]0,L[))$, $T(u(\cdot)) \in L^1_{loc, w}(0, T, BV(]0,L[))$ for all $T \in \mathcal
T_r$, and $\z(t):= \a(u(t), \partial_x u(t)) \in L^1(Q_T)$,  such that:

\begin{itemize}
\item[(i)] the
time derivative $u_t$ of $u$ in $\mathcal{D}'(Q_T)$ belongs to
$[L^1(0, T; BV(]0,L[))]^*$ and satisfies
\begin{equation}\label{IBPFTTT}
\int_0^T \langle u_t(t), \psi(t) \rangle \, dt  = -\int_0^T
\int_0^L u(t,x) \Theta(t,x) \, dxdt
\end{equation}
for all test function $\psi \in L^1(0, T; BV(]0,L[))$ compactly supported in time such that
$\psi(t) = \int_0^t \Theta(s) \, ds$ and  $\Theta \in
L_w^1(0,T;BV(]0,L[)) \cap L^{\infty}(Q_T)$.
\item[(ii)] \  $D_x \z =u_t$ in $\mathcal{D}'(Q_T)$, and for any
$w \in L^1(0, T; BV(]0,L[))$, the distribution $\z D_xw$ defined by
(\ref{DefMRAd}) is a Radon measure in $Q_T$ and verifies, for all $w \in L^1(0, T; BV(]0,L[))$, the following
integration by parts formula
\begin{equation}\label{IBPF}
   \int_{Q_T} \z D_xw +  \langle  u_t, w \rangle=
    \beta \int_0^T  w(t, 0_+) \, dt- c \int_0^T u(t, L_-) w(t, L_-) \, dt.
\end{equation}
\item[(iii)] \  the following inequality is satisfied
$$
\displaystyle \int_{Q_T} \eta
h_S(u,DT(u)) \, dt + \int_{Q_T} \eta h_T(u,DS(u)) \, dt
 \leq
 $$
 $$
 \displaystyle\int_{Q_T}J_{TS}(u) \partial_t\eta\ dxdt -
\int_{Q_T} \a(u, \partial_x u) \cdot \partial_x \eta \ T(u)
S(u) \ dxdt
$$
 for truncatures $S , \, T \in {\mathcal T}^+$  and any $\eta \in C^\infty (Q_T)$ of
 compact support.
\end{itemize}
 }
\end{definition}

In the following result we get a positive lower bound for $u(t,
0_+)$.

\begin{lemma}\label{mayorcero} If $u$ is an entropy solution of (\ref{Dirichletproblem}) in
$Q_T = (0,T)\times ]0,L[$, then
\begin{equation}\label{Erty333} u(t, 0_+) \geq \frac{\beta}{c}> 0, \ \ \ \
 \qquad \hbox{for almost all } t \in ]0,T[.
\end{equation}
\end{lemma}
\noindent{\bf Proof.} For any $n \in \N$, let
$v_n$ be the function defined by zero in $]1/n,L]$, $1$ at $x=0$, and a straight line joining both values at the rest of the points.
Being $0 \leq \phi \in {\mathcal D}(]0, T[)$ fixed and taking $w$ in
(\ref{IBPF}) as $w_n(t):= \phi(t) v_n$, we get
\begin{equation}\label{E2mayorcero}
   \int_{Q_T} \z Dw_n +  \langle  u_t, w_n \rangle =
    \beta \int_0^T  \phi(t) \, dt.
\end{equation}
By (\ref{IBPFTTT}), we have
$$\langle  u_t, w_n\rangle = - \int_0^T
\phi'(t) \int_0^L u(t,x) v_n(x) \, dx dt,$$ so by the Dominate
Convergence Theorem,
\begin{equation}\label{E3mayorcero}
\lim_{n \to \infty}  \langle  u_t, w_n \rangle  = 0.
\end{equation}
On the other hand, given $\varphi \in  {\mathcal D}(Q_T)$, we have
$$\langle \z D_xw_n, \varphi \rangle  =
\int_0^T  \phi(t)  \int_0^L \z(t,x)
    \varphi(t,x) v'_n(x) \, dx dt.
    $$
    Hence,
    \begin{equation}\label{E4mayorcero}
   \int_{Q_T} \z(t,x) D_xw_n(t,x)  = - \int_0^T n \phi(t)
   \int_0^{\frac{1}{n}} \z(t,x) \, dx dt.
\end{equation}
Now, by (\ref{E2mayorcero}), (\ref{E3mayorcero}) and
(\ref{E4mayorcero}), we get
$$ \beta \int_0^T  \phi(t) \, dt = - \lim_{n \to \infty} \int_0^T \phi(t)  n
   \int_0^{\frac{1}{n}} \z(t,x) \, dx dt.$$
Then, since $\vert \z(t,x) \vert \leq c u(t,x)$, by Fatou's Lemma
we obtain that

$$\beta \int_0^T  \phi(t) \, dt \leq c  \int_0^T \phi(t) \left[\lim_{n \to \infty} n
   \int_0^{\frac{1}{n}} u(t,x) \, dx  \right]dt = c  \int_0^T \phi(t) u(t,0_+) \, dt$$
from where it follows (\ref{Erty333}). \hfill$\Box$

\begin{remark}\label{timereg}{\rm Let $u$ a
 bounded entropy solution of (\ref{Dirichletproblem}) in $Q_T$. In the proof of the next
 result we need the following time regularization.
  For that, given $\phi\in \mathcal D(]0,T[)$ and  $w \in
L_{loc}^1(0,T; BV(]0,L[))$, we define $(\phi w)^{\tau}$, as the Dunford integral (see
\cite{DU})
\begin{equation}\nonumber\label{SuperFuensanta}(\phi w)^{\tau}(t):=
\frac{1}{\tau} \int_{t - \tau}^t \phi (s) w(s) \ ds \in BV(]0,L[)^{**},
\end{equation}  that is
$$\langle (\phi w)^{\tau}(t), \eta \rangle =  \frac{1}{\tau} \int_{t - \tau}^t \langle \phi (s) w(s), \eta \rangle \
ds \ \ \ \forall \, \eta \in BV(]0,L[)^*.$$ In \cite{ACM4:01} it is shown that $(\phi
w)^{\tau} \in C([0,T]; BV(]0,L[))$. If $u$ is an entropy solution
of (\ref{Dirichletproblem}) and $p \in {\mathcal T}^+$, it is easy to see that
$$
|D_x(\phi  p(u))^\tau(t)|(]0,L[) \leq \frac{1}{\tau}\int_{t-\tau}^t
|D_x(\phi(s) p(u(s)))|(]0,L[)\, ds.
$$
Then, by the lower-semi-continuity of the total
variation respect to the $L^1$-convergence, we have
$$
|D_x(\phi(t) p(u(t)))|(]0,L[) \leq \liminf_{\tau \to 0} |D_x(\phi p(u))^\tau(t)|(]0,L[)
$$
$$
\leq \limsup_{\tau \to 0}\frac{1}{\tau}\int_{t-\tau}^t
|D_x(\phi(s) p(u(s))|(]0,L[)\, ds.
$$
 Since the map $t \mapsto |D_x(\phi(t)  p(u(t)))|(]0,L[)$
belongs to  $L^1_{\rm loc}([0,T])$, we have that almost all $t\in [0,T]$ is a Lebesgue
point of this map. So, for almost all $t\in [0,T]$, we have
$$
\frac{1}{\tau}\int_{t-\tau}^t
|D_x(\phi(s) p(u(s))|(]0,L[)\, ds\stackrel{\tau\to 0}\longrightarrow \vert D_x(\phi(t) p(u(t))
\vert(]0,L[),
$$
and consequently,

\begin{equation}\nonumber\label{conmess}|D_x(\phi p(u))^\tau(t)|(]0,L[) \stackrel{\tau\to 0}\longrightarrow \vert
D_x(\phi(t) p(u(t)) \vert(]0,L[) \ a.e.\ t.\end{equation}

}
\end{remark}

Respect to the existence and uniqueness of bounded entropy solutions we have the
following result.

\begin{theorem}\label{EXistparabolic}
 For any initial datum
  $0 \leq u_0 \in   L^{\infty}(]0,L[)$ there exists a
 unique bounded entropy solution $u$ of
(\ref{Dirichletproblem}) in $Q_T = ]0,T[\times ]0,L[$ for every $T
> 0$ such that $u(0) = u_0$. Moreover, if $u(t)$, $\overline{u}(t)$ are bounded entropy solutions of
(\ref{Dirichletproblem}) in $Q_T = ]0,T[\times ]0,L[$ co\-rres\-ponding to initial data
$u_0$, $\overline{u}_0 \in L^{\infty}(]0,L[)^+$, respectively, then
\begin{equation}\nonumber
\label{CPentropys}\Vert (u(t) - \overline{u}(t))^+ \Vert_1 \leq \Vert (u_0 -
\overline{u}_0)^+ \Vert_1 \ \ \ \ \ \ {\rm for \ all} \ \ t \geq 0.
\end{equation}
In particular, we have uniqueness of bounded entropy solutions of
(\ref{Dirichletproblem}).
\end{theorem}
\noindent{\bf Proof.} \noindent{\it The comparison principle.} Let $b > a > 2 \epsilon
> 0$, $T(r):=T_{a,b}(r)-a$. We need to consider truncature functions of the form
$S_{\epsilon, l}(r):= T_{\epsilon}(r - l)^+ = T_{l, l + \epsilon}(r) - l \in {\mathcal
T}^+$,
 and $S_\epsilon^l(r):= T_{\epsilon}(r - l)^- + \epsilon =
T_{l- \epsilon, l}(r) + \epsilon - l \in {\mathcal T}^+,$ where $l \geq 0$. Observe that
$ S_\epsilon^l(r) = - T_\epsilon(l-r)^+ + \epsilon. $  Let us denote
$$
J^+_{T,\epsilon,l}(r) = \int_0^r T(s) T_{\epsilon}(s-l)^+ \, ds, $$ $$
J^-_{T,\epsilon,l}(r) = \int_0^r T(s) T_{\epsilon}(s-l)^- \, ds = -\int_0^r T(s)
T_{\epsilon}(l-s)^+ \, ds.
$$
Then, $J_{TS_{\epsilon, l}}(r) = J^+_{T,\epsilon,l}(r)$ and $J_{TS_\epsilon^l}(r) =
J^-_{T,\epsilon,l}(r)+ \epsilon J_T(r).$

Let $u$, $\overline{u}$ be two entropy solutions of (\ref{Dirichletproblem})
corresponding to the initial conditions $u_0, \overline{u_0} \in
\left(L^{1}(]0,L[)\right)^+$, respectively. Then,  if $\z(t) := \a(u(t),
\partial_x u(t))$, $\zbar (t) := \a(\ubar(t),
 \partial_x \ubar(t))$, and   $l_1, l_2  > \epsilon$, we have
\begin{equation}\label{UE3}
\begin{array}{ll}
- \displaystyle \int_0^T \int_0^L J^+_{T,\epsilon,l_1}(u(t)) \partial_t\eta(t) \ dxdt \\ \\+
\displaystyle \int_0^T \int_0^L \eta(t) [h_T(u(t),
D_xS_{\epsilon,l_1}(u(t))) + h_{S_{\epsilon,l_1}}(u(t),D_xT(u(t)))] \ dt \\ \\
+ \displaystyle \int_0^T \int_0^L \z(t)   \partial_x \eta(t) \ T(u(t))
S_{\epsilon,l_1}(u(t)) \ dxdt
 \leq 0,
\end{array}
\end{equation}
and
\begin{equation}\label{UE4}
\begin{array}{l}
- \displaystyle \int_0^T \int_0^L J^-_{T,\epsilon,l_2}( \overline{u}(t)) \partial_t \eta \ dxdt -
\displaystyle \epsilon \int_0^T
\int_0^L J_{T}( \overline{u}(t)) \partial_t \eta(t) \ dxdt \\
\\+ \displaystyle \int_0^T \int_0^L \eta(t)
[h_T(\overline{u}(t), D_xS^{l_2}_{\epsilon} (\overline{u}(t)) )
+  h_{S^{l_2}_{\epsilon}}(\overline{u}(t),D_xT(\overline{u}(t)))] \ dt \\
\\+
\displaystyle \int_0^T \int_0^L \overline{\z}(t)   \partial_x \eta(t) \
T(\overline{u}(t)) S^{l_2}_{\epsilon}(\overline{u}(t))
 \ dxdt \leq 0,
\end{array}
\end{equation}
for all $\eta \in C^{\infty} (Q_T)$, with \ $\eta \geq 0$, $\eta(t, x) = \phi(t)
\rho(x)$, being \ $\phi \in {\mathcal D}(]0, T[)$, \ $\rho \in {\mathcal D} (]0,L[)$.

We choose two different pairs of variables $(t,x)$, $(s, y)$ and consider $u$, $\z$ as
functions in $(t,x)$, $\overline{u}$, $\overline{\z}$ in $(s, y)$. Let $0 \leq \phi \in
{\mathcal D}(]0, T[)$, $\psi\in \mathcal D(]0,L[)$,  $\rho_m$ and $\tilde{\rho}_n$
sequences of mollifier in $\R$. Define
$$\eta_{m,n}(t, x, s, y):= \rho_m(x - y) \tilde{\rho}_n(t - s) \phi
\left(\frac{t+s}{2} \right)\psi\left(\frac{x+y}{2}\right).$$ For $(s,y)$ fixed, if we
take in (\ref{UE3}) $l_1 = \overline{u}(s, y)$, we get
\begin{equation}\label{UE5}
\begin{array}{llll}
- \displaystyle \int_0^T \int_0^L
J^+_{T,\epsilon,\overline{u}(s, y)}(u(t,x)) \partial_t \eta_{m,n} \ dx dt  \\ \\
+ \displaystyle \int_0^T \int_0^L \eta_{m,n} [h_T(u(t,x), D_xS_{\epsilon,\overline{u}(s,
y)}(u(t,x)) ) + h_{S_{\epsilon, \overline{u}(s, y)}}(u(t,x),D_xT(u(t,x)))]
\ dt  \\
\\ + \displaystyle \int_0^T \int_0^L \z(t,x)
 \partial_x \eta_{m,n} \ T(u(t,x)) \ S_{\epsilon,
\overline{u}(s, y)}(u(t,x)) \ dx dt \leq 0.
\end{array}
\end{equation}
Similarly, for $(t,x)$ fixed, if we take in (\ref{UE4}) $l_2 = u(t, x)$  we get
\begin{equation}\label{UE6}
\begin{array}{llll}
- \displaystyle \int_0^T \int_0^L J^-_{T,\epsilon,u(t,x)}( \overline{u}(s, y))
\partial_s \eta_{m,n} \ dy ds  - \displaystyle \epsilon \int_0^T \int_0^L J_{T}(\overline{u}(s,
y)) \partial_s \eta_{m,n} \ dy ds
\\ \\
+ \displaystyle \int_0^T \int_0^L \eta_{m,n} [h_T(\overline{u}(s, y), D_y
S^{u(t,x)}_{\epsilon}( \overline{u}(s, y)) ) +
h_{S^{u(t,x)}_{\epsilon}}(\overline{u}(s,y),D_yT(\overline{u}(s,y)))]
\,  ds  \\
\\ + \displaystyle \int_0^T \int_0^L \overline{\z}(s, y)
\partial_y \eta_{m,n} \ T(\overline{u}(s,y)) \
S^{u(t,x)}_{\epsilon}(\overline{u}(s, y)) \, dy ds   \leq 0.
\end{array}
\end{equation}

We integrate (\ref{UE5}) in $(s, y)$, (\ref{UE6}) in $(t, x)$, and add the two
inequalities. Using that $a > 2 \epsilon$, and since

$$
\int_{Q_T \times Q_T} \eta_{m,n} h_{S_{\epsilon,\overline{u}(s,
y)}}(u(t,x),D_xT(u(t,x)))\ dsdtdy \geq 0$$ and
$$\displaystyle \int_{Q_T \times Q_T}
\eta_{m,n} h_{S^{u(t,x)}_{\epsilon}}(\overline{u}(s,y),D_yT(\overline{u}(s,y)))\ dsdtdx
\geq 0,
$$
we get
\begin{equation}\label{UE8}
\begin{array}{l}
-  \displaystyle\int_{Q_T \times Q_T} \left(J^+_{T,\epsilon,\overline{u}(s, y)}(u(t,x))
\partial_t \eta_{m,n} + J^-_{T,\epsilon,u(t,x)}( \overline{u}(s, y)) \partial_s \eta_{m,n} \right)\
dsdtdydx
\\ \\
 -
\epsilon \displaystyle\int_{Q_T \times Q_T} J_{T}(\overline{u}(s, y)) \partial_s \eta_{m,n} \ dsdtdydx
 \\ \\
 +
 \displaystyle\int_{Q_T \times Q_T} \eta_{m,n} h_T(u(t,x), D_x S_{\epsilon, \overline{u}(s,
y)}(u(t,x)))\ dsdtdy
  \\ \\ +
\displaystyle \int_{Q_T \times Q_T} \eta_{m,n} h_T(\overline{u}(s, y), D_y S^{u(t,x)}_{\epsilon}(\overline{u}(s,
y)))\ dsdtdx
\\ \\
- \displaystyle \int_{Q_T \times Q_T} \overline{\z}(s, y)
 \partial_x \eta_{m,n}  T(\overline{u}(s, y)) S^{u(t,x)}_{\epsilon}(\overline{u}(s, y)) \ dsdtdydx
\\ \\ -  \displaystyle\int_{Q_T \times Q_T} \z(t, x)
\partial_y \eta_{m,n} T(u(t,x)) S_{\epsilon, \overline{u}(s,
y)}(u(t,x) ) \ dsdtdydx
\\ \\
+ \displaystyle \int_{Q_T \times Q_T}T^+_\epsilon (u(t,x)-\overline
u(s,y))[T(u(t,x))\z(t,x)-T(\overline u(s,y))\overline\z(s,y)]
\\ \\ \times (\partial_x \eta_{m,n}+\partial_y\eta_{m,n}) \ dsdtdydx
\\ \\
\displaystyle +\epsilon\int_{Q_T \times Q_T} T(\overline
u(s,y))\overline\z(s,y)(\partial_x \eta_{m,n}+\partial_y\eta_{m,n}) \ dsdtdydx  \leq 0.
\end{array}
\end{equation}
Let $I_2$ be the sum of the third up to the sixth terms of the above inequality. From now
on, since $u$, $\z$ are always functions of $(t, x)$, and $\overline{u}$, $\overline{\z}$
are always functions of $(s,y)$, to make our expression shorter, we shall omit the
arguments except when they appear as sub-index and in some additional cases where we find
it useful to remind them. We also omit the differentials of the integrals.

Working as in the proof of uniqueness  of Theorem 3 in \cite{ARMA}, we obtain that $
\frac{1}{\epsilon}  I_2\geq \Vert \phi\Vert_\infty \Vert \psi\Vert_\infty o(\epsilon). $
Hence, by (\ref{UE8}), it follows that
$$- \displaystyle \int_{Q_T \times Q_T}
\left(J^+_{T,\epsilon,\overline{u}}(u) \partial_t \eta_{m,n} + J^-_{T,\epsilon, u}(
\overline{u}) \partial_s \eta_{m,n} \right)$$
$$ +\int_{Q_T \times Q_T}T^+_\epsilon (u-\overline
u)[T(u)\z-T(\overline u)\overline\z](\partial_x \eta_{m,n}+\partial_y\eta_{m,n})
$$
$$+\epsilon\int_{Q_T \times Q_T} T(\overline
u)\overline\z(\partial_x \eta_{m,n}+\partial_y\eta_{m,n}) \leq \epsilon
o(\epsilon)+\epsilon \int_{Q_T \times Q_T} J_{T}(\overline{u}) \partial_s\eta_{m,n}.$$
Then, dividing by $\epsilon$ and letting $\epsilon\to 0$ we get
$$-\int_{Q_T \times Q_T} \left(J^+_{T,{\rm sign},\overline
u}(u)\partial_t \eta_{m,n}+J^-_{T,{\rm sign},u}(\overline u)\partial_s \eta_{m,n}\right)$$$$
+\int_{Q_T \times Q_T}{\rm sign}_0^+ (u-\overline u)[T(u)\z-T(\overline u)\overline\z](\partial_x
\eta_{m,n}+\partial_y\eta_{m,n})
$$$$+\int_{Q_T \times Q_T}T(\overline
u)\overline\z(\partial_x \eta_{m,n}+\partial_y\eta_{m,n})\leq \int_{Q_T \times Q_T}
J_{T}(\overline{u}) \partial_s \eta_{m,n}
$$
  where
$$
J^+_{T,sign,l}(r) = \int_0^r T(r') sign^+_0(r'-l) dr' \quad \hbox{$l\in\R$, $r\geq 0$}
$$
and
$$
J^-_{T,sign,l}(r) = \int_0^r T(r') sign^-_0(r'-l) dr' \quad \hbox{$l\in\R$, $r\geq 0$.}
$$
Now, letting $m\to \infty$, we obtain
$$-\int_0^T \int_0^T \int_0^L\left(J^+_{T,{\rm sign},\overline
u(s,x)}(u(t,x))\partial_t \1_{n}+J^-_{T,{\rm sign},u(t,x)}(\overline
u(s,x))\partial_s \1_{n}\right)
$$
$$
 +\int_0^T \int_0^T \int_0^L{\rm
sign}_0^+ (u(t,x)-\overline u(s,x))[T(u(t,x))\z(t,x)-T(\overline u(s,x))\overline\z(s,x)]
\partial_x \1_{n}
$$
$$
+\int_0^T \int_0^T \int_0^L T(\overline
u(s,x))\overline\z(s,x)\partial_x\1_{n} \leq \int_0^T \int_0^T \int_0^LJ_{T}(\overline
u(s,x))\partial_s \1_{n}
$$
 where $\1_n(t,s,x):=\tilde\rho_n(t-s)\phi(\frac{t+s}{2})\psi(x)$.
We set $\psi=\psi_k\in\mathcal D(]0,L[)\uparrow \1_{]0,L[}$ in the last expression and
taking limit as $k \to + \infty$, we have
\begin{equation}
\label{psim} \begin{array}{l}-\displaystyle\int_0^T \int_0^T \int_0^L \left(J^+_{T, {\rm
sign},\overline u(s,x)}(u(t,x)) \partial_t\kappa_n(t,s)+J^-_{T, {\rm sign},u(t,x)}(\overline
u(s,x))\partial_s \kappa_n(t,s)\right)
\\ \\ \displaystyle+\lim_{k\to
+\infty}\int_0^T \int_0^T  \int_0^L \kappa_n(t,s){\rm sign}_0^+
(u(t,x)-\overline u(s,x))T(u(t,x))\z(t,x)\partial_x \psi_k(x)  \\ \\
\displaystyle-\lim_{k\to +\infty}\int_0^T \int_0^T  \int_0^L \kappa_n(t,s){\rm sign}_0^+
(u(t,x)-\overline u(s,x))T(\overline
u(s,x))\overline\z(s,x))\partial_x \psi_k(x) \\ \\
\displaystyle+\lim_{k\to +\infty}\int_0^T \int_0^T \int_0^L\kappa_n(t,s)T(\overline
u(s,x)\overline\z(s,x))\partial_x \psi_k(x) \\ \\ \displaystyle\leq \int_0^T \int_0^T
\int_0^L J_{T}(\overline u(s,x))\partial_s \kappa_n(t,s), \end{array}
\end{equation}
where $\kappa_n(t,s):=\tilde\rho_n(t-s)\phi(\frac{t+s}{2})$.

Let us study the second and the third term of the above expression. Let
$$I_k:= \int_0^T \int_0^T  \int_0^L \kappa_n(t,s){\rm sign}_0^+
(u(t,x)-\overline u(s,x))T(u(t,x))\z(t,x)\partial_x \psi_k(x)
$$
$$
= \int_0^T \int_0^T  \int_0^L \kappa_n(t,s){\rm sign}_0^+
(u(t,x)-\overline u(s,x))T(u(t,x))\z(t,x)\partial_x (\psi_k(x)- 1).
$$
Let $H_n(s,r):=
\kappa_n(r,s){\rm sign}_0^+ (u(r)-\overline u(s))T(u(r))$. For $\tau > 0$, we define the
function $(\kappa_n(s))^{\tau}$, as the Dunford integral (see Remark \ref{timereg})
$$
(\kappa_n(s))^{\tau}(t):= \frac{1}{\tau}\int_t^{t+\tau}H_n(s,r) \ dr.
$$
 Then,
$$
I_k = \lim_{\tau \to 0} \int_0^T \int_0^T  \int_0^L (\kappa_n(s))^{\tau}(t)\z(t,x)
\partial_x [\psi_k(x)- 1] \, dxdtds
$$
$$
= - \lim_{\tau \to 0} \int_0^T \int_0^T \int_0^L [\psi_k(x)- 1] \z(t,x) D_x((\kappa_n(s))^{\tau}(t)))
\, ds dt
$$
$$
 - \lim_{\tau \to 0} \int_0^T   \langle u_t, (\kappa_n(s))^{\tau}(\psi_k(x)- 1) \rangle
\, ds$$
$$
+ c \lim_{\tau \to 0}  \int_0^T  \int_0^T
 u(t, L_{-})(\kappa_n(s))^{\tau}(t)(L_{-}) \, dt ds
$$
  $$
   - \beta \lim_{\tau \to 0}  \int_0^T  \int_0^T
 (\kappa_n(s))^{\tau}(t)(0_+) \, dt ds
 = I^1_k  + I^2_k + I^3_k + I^4_k.
$$
Notice that
$$I^3_k = c \int_0^T  \int_0^T
 u(t, L_{-}) \kappa_n(t,s){\rm sign}_0^+ (u(t,L_{-} ) -\overline u(s, L_{-}))T(u(t, L_{-}) )
\,  dt ds$$ and
$$I^4_k = - \beta \int_0^T  \int_0^T
 \kappa_n(t,s){\rm sign}_0^+ (u(t, 0_+) -\overline u(s, 0_+))T(u(t, 0_+) )
\,  dt ds.$$
By Remark \ref{timereg},  we get
\begin{equation}\label{convvRtyDD}
|D_x((\kappa_n(s))^{\tau}(t))|(]0,L[) \stackrel{\tau\to 0}\longrightarrow \vert
D_x(\kappa_n(t,s){\rm sign}_0^+ (u(t)-\overline u(s))T(u(t))) \vert(]0,L[).
\end{equation}
Using (\ref{convvRtyDD}), we get

$$\vert I^1_k \vert \le c
\Vert u \Vert_{L^{\infty}(Q_T)}\int_0^T \int_0^T \int_0^L (1 - \psi_k(x)) \vert D_x(\kappa_n(t,s){\rm
sign}_0^+ (u(t)-\overline u(s))T(u(t))) \vert \ dtds,$$ which implies $\lim_{k \to
\infty} I^1_k = 0.$ Let us deal with $I^2_k$. We have

$$I^2_k  = \lim_{\tau \to 0}  \int_0^T  \int_0^T \int_0^L u(t,x) \frac{H_n(s,t+\tau)- H_n(s,t)}{\tau}
(\psi_k(x) - 1) \, dx dt ds.$$ Let
$$q(\tau):= {\rm sign}_0^+(\tau - \overline u(s,x))T(\tau), \ \ \
Q(r):= \int_0^r q(\tau) \, d\tau.$$ Since $q$ is non-decreasing, $Q(r) - Q(\overline{r})
\leq q(r)(r -\overline{r}).$ Then, changing variables, since $H_n(s,t) = q(u(t))
\kappa_n(t,s)$

\begin{equation}
\label{truco_estandar}
\begin{array}{l}
I^2_k = \displaystyle\lim_{\tau \to 0}  \int_0^T  \int_0^T \int_0^L (1 - \psi_k(x)) \frac{u(t,x) - u(t
- \tau,x)}{\tau} H_n(s,t) \, dx dt ds
\\ \\
 \geq \displaystyle \lim_{\tau \to 0} \int_0^T \int_0^T \int_0^L (1 -
\psi_k(x))
 \kappa_n(t,s) \frac{Q(u(t,x)) - Q(u(t-\tau,x))}{\tau} \, dx dt ds
 \\ \\
 =\displaystyle \lim_{\tau \to 0}  \int_0^T  \int_0^T \int_0^L (1 -
\psi_k(x)) Q(u(t,x)) \frac{\kappa_n(t,s) - \kappa_n(t+ \tau,s)}{\tau} \, dx dt ds
\\ \\
 = - \displaystyle\int_0^T  \int_0^T \int_0^L (1 - \psi_k(x)) Q(u(t,x))
\partial_t \kappa_n(t,s)  \, dx dt ds,
\end{array}
\end{equation}
from where it follows that $\lim_{k \to \infty} I^2_k \geq  0.$ Taking into account the
above facts, we get
\begin{equation}\label{Santica1}\begin{array}{l}
\displaystyle \lim_{k \to \infty} I_k \geq - \beta \displaystyle \int_0^T  \int_0^T
 \kappa_n(t,s){\rm sign}_0^+ (u(t, 0_+) -\overline u(s, 0_+))T(u(t, 0_+) )
\,  dt ds \\ \\ + c \displaystyle\int_0^T  \int_0^T
 u(t, L_{-}) \kappa_n(t,s){\rm sign}_0^+ (u(t,L_{-} ) -\overline u(s, L_{-}))T(u(t, L_{-}) )
\,  dt ds.
\end{array}
\end{equation}
 Working similarly, we obtain
\begin{equation}\label{Santica2} \begin{array}{l} -\displaystyle
\lim_{k \to \infty} \int_0^T  \int_0^T \int_0^L\kappa_n(t,s){\rm sign}_0^+ (u(t,x)- \overline
u(s,x))T(\overline u(s,x))\overline\z(s,x) \partial_x \psi_k(x) \\ \\
\geq  \beta \displaystyle\int_0^T  \int_0^T
  \kappa_n(t,s){\rm
sign}_0^+ (u(t,0_+)-\overline u(s,0_+))T(\overline{u}(s,0_+)) \, dt ds \\ \\
\displaystyle -c \int_0^T  \int_0^T
  \overline{u}(s, L_{-})\kappa_n(t,s){\rm
sign}_0^+ (u(t,L_{-})-\overline u(s,L_{-}))T(\overline{u}(s,L_{-})) \, dt ds.
\end{array}
\end{equation}
Analogously,
\begin{equation}\label{Santica3}\begin{array}{l}\displaystyle
\lim_{k\to +\infty}\int_0^T  \int_0^T \int_0^L \kappa_n(t,s)T(\overline
u(s,x))\overline\z(s,x)\partial_x \psi_{k}(x) \\ \\ \displaystyle\geq  c \int_0^T \int_0^T
\overline{u}(s,L_{-}) \kappa_n(t,s)T(\overline{u}(s, L_-)) \,  dt ds- \beta \int_0^T \int_0^T
\kappa_n(t,s)T(\overline{u}(s,0_+)) \,  dt ds.\end{array}
\end{equation}

From (\ref{psim}), by (\ref{Santica1}), (\ref{Santica2}) and (\ref{Santica3}), we have
\begin{equation}
\label{psimTTT}
\begin{array}{l}
-\displaystyle\int_0^T  \int_0^T \int_0^L \left(J^+_{T,{\rm sign},\overline
u(s,x)}(u(t,x))\partial_t \kappa_n(t,s)+J^-_{T,{\rm sign},u(t,x)}(\overline
u(s,x))\partial_s \kappa_n(t,s)\right)\ dtdsdx \\ \\ + c \displaystyle\int_0^T  \int_0^T
 u(t, L_{-}) \kappa_n(t,s){\rm sign}_0^+ (u(t,L_{-} ) -\overline u(s, L_{-}))T(u(t, L_{-}) )
\,  dt ds \\ \\
\displaystyle -c \int_0^T  \int_0^T
  \overline{u}(s, L_{-})\kappa_n(t,s){\rm
sign}_0^+ (u(t,L_{-})-\overline u(s,L_{-}))T(\overline{u}(s,L_{-})) \, dt ds
\\ \\
 \displaystyle - \beta \int_0^T  \int_0^T
 \kappa_n(t,s){\rm sign}_0^+ (u(t, 0_+) -\overline u(s, 0_+))T(u(t, 0_+)
 )
\,  dt ds
\\ \\
+\displaystyle \beta \int_0^T  \int_0^T
  \kappa_n(t,s){\rm
sign}_0^+ (u(t,0_+)-\overline u(s,0_+))T(\overline{u}(s,0_+)) \,
dt ds \\ \\
\displaystyle + c \int_0^T \int_0^T \overline{u}(s,L_{-})
\kappa_n(t,s)T(\overline{u}(s,L_{-})) \, dt ds- \beta \int_0^T \int_0^T
\kappa_n(t,s)T(\overline{u}(s,0_+)) \,  dt ds
 \\ \\ \displaystyle\leq \displaystyle \int_0^T  \int_0^T \int_0^L J_{T}(\overline u(s,x))\partial_s \kappa_n(t,s)\
dtdsdx.
\end{array}
\end{equation}
By Lemma \ref{mayorcero}, we have
\begin{equation}\label{Erty22} u(t, 0_+) \geq \frac{\beta}{c}> 0, \ \ \ \ \overline
u(s,0_+)\geq \frac{\beta}{c}
> 0 \qquad \hbox{for almost every} t, s > 0.
\end{equation}
 Letting $a\to 0$, dividing by $b$ and letting $b\to 0$ in (\ref{psimTTT}), we obtain,
$$-\int_0^T  \int_0^T
\int_0^L (u(t,x)-\overline u(s,x))^+(\partial_t\kappa_{n}(t,s)+\partial_s\kappa_{n}(t,s)) \ dtdsdx$$
$$ + c \displaystyle\int_0^T  \int_0^T
 u(t, L_{-}) \kappa_n(t,s){\rm sign}_0^+ (u(t,L_{-} ) -\overline u(s, L_{-})){\rm sign}_0^+ (u(t, L_{-}) )
\,  dt ds $$ $$ -c \int_0^T  \int_0^T
  \overline{u}(s, L_{-})\kappa_n(t,s){\rm
sign}_0^+ (u(t,L_{-})-\overline u(s,L_{-})){\rm sign}_0^+(\overline{u}(s,L_{-})) \, dt ds
$$
$$ - \beta \int_0^T  \int_0^T
  \kappa_n(t,s){\rm sign}_0^+ (u(t, 0_+)-\overline
u(s, 0_+)) [ {\rm sign}_0(u(t, 0_+) - {\rm sign}_0(\overline{u}(s,0_+))]  dt ds
$$
$$+ c \int_0^T \int_0^T \overline{u}(s,L_{-}) \kappa_n(t,s){\rm sign}_0^+(\overline{u}(s,L_{-}))
\, dt ds- \beta \int_0^T \int_0^T \kappa_n(t,s)  \,  dt ds  $$ $$\leq \int_0^T  \int_0^T
\int_0^L \overline u(s,x)\partial_s\kappa_{n}(t,s)\ dtdsdx.$$ Having in mind (\ref{Erty22}),
the fourth term of the above expression vanishes. Moreover, the sum of the second and
third term is  non-negative. On the other hand, since $ \overline{u}_s =  D_x \,
(\overline{\z})$ in the sense given in (ii) of Definition \ref{DefESevP},
 $$\int_0^T  \int_0^T
\int_0^L\overline u(s,x)\partial_s\kappa_{n}(t,s)  \,  dxdt ds= - \int_0^T \langle  \overline{u}_s,
\kappa_n(\cdot,t) \rangle \,  dt $$ $$= c \int_0^T \int_0^T \overline{u}(s,L_{-})
\kappa_n(t,s) \, dt ds- \beta \int_0^T \int_0^T \kappa_n(t,s)  \, dt ds.$$ Therefore,
$$-\int_0^T  \int_0^T \int_0^L(u(t,x)-\overline
u(s,x))^+(\partial_t\kappa_{n}(t,s)+\partial_s\kappa_{n}(t,s))\ dtdsdx \leq 0.$$ Letting $n\to\infty$, $$-
\int_0^T \int_0^L (u(t,x)-\overline u(t,x))^+\phi'(t)\ dxdt\leq 0.$$ Since this is true
for all $0 \leq \phi \in {\mathcal D}(]0, T[)$, we have
$$\frac{d}{dt} \int_0^L
(u(t,x) - \overline{u}(t, x))^+ \, dx \leq 0.$$ Hence
$$\int_0^L (u(t,x) - \overline{u}(t, x) )^+ \, dx
\leq \int_0^L (u_0(x) - \overline{u}_0(x))^+ \, dx \ \ \ \ \ \ {\rm for \ all} \ \ t \geq
0,$$ which finishes the uniqueness part.

\medskip

 \noindent{\it Existence of bounded entropy solution.} Given $0 \leq u_0 \in L^{1}(]0,L[)$, let $u(t) = T_{\beta}(t) u_0$, being
$(T_{\beta}(t))_{t \geq 0}$ the semigroup in $L^1(]0,L[)^+$ generated by the accretive
operator ${\mathcal B}_{\beta}$. Then, according to the general theory of nonlinear
semigroups (\cite{BCrP}), we have that $u(t)$ is a mild-solution of the abstract Cauchy
problem
$$
 u^{\prime} (t) + \mathcal{B}_{\beta}u(t) \ni 0, \ \ \ \ u(0) = u_0.
$$

Let us prove that, assuming $0 \leq u_0 \in L^{\infty}(]0,L[)$, then $u$ is a bounded
entropy solution of (\ref{Dirichletproblem}) in $Q_T$. We divide the proof of existence
in several steps.

\noindent{\it Step 1. Approximation with Crandall-Ligget's
scheme.} Let $T>0$, $K \geq 1$, $\Delta t = \frac{T}{K}$, $t_n =
n\Delta t$, $n=0,\ldots,K$. We define inductively $u^{n+1}$,
$n=0,\ldots,K-1$  to be the unique entropy solution of
\begin{equation}
\label{CL1} \left\{\begin{array}{ll} \displaystyle
\frac{u^{n+1}-u^n}{\Delta t} - \left(\a(u^{n+1},(u^{n+1})')
\right)' = 0 &\hspace{0.1cm}\hbox{in \hspace{0.1cm} $]0, L[$}
\\ \\ -  \a(u^{n+1}(0),(u^{n+1})'(0)) = \beta > 0 \ \ \hbox{and} \ \
u^{n+1}(L_{-})=0,
\end{array} \right.
\end{equation}
where $u^0 = u_0$.

 If we define
$$u^K(t) := u^0 \1_{[0,t_1]}(t)+\sum_{n=1}^{K-1} u^n
\1_{]t_n,t_{n+1}]}(t),$$ by Crandall-Liggett's Theorem, we get that
$u^K$ converges uniformly to $ u \in C([0,T],L^1(]0, L[))$, as $ K \to
\infty$.

We also define
$$\xi^K(t):= \sum_{n=0}^{K-1}
\frac{u^{n+1}-u^n}{\Delta t} \1_{]t_n,t_{n+1}]}(t)$$ and
$$\z^K(t) := \a(u^1, (u^1)')
\1_{[0,t_1]}(t)+\sum_{n=1}^{K-1} \a(u^{n+1}, (u^{n+1})^\prime)
\1_{]t_n,t_{n+1}]}(t).$$ Since $u^{n+1}$ is the entropy solution
of (\ref{CL1}),  we have
\begin{equation}\label{E1exist}
\xi^K(t) = D_x \z^K(t) \ \ {\rm in} \ \ {\mathcal D}^{\prime}(]0, L[), \ \ \forall t \in
]0,T]
\end{equation}
\begin{equation}
\label{condel2existt} \z^K(t)(L) = - c u^K(t+ \Delta t)(L_{-}), \ \ \forall t \in  ]0, T
-\step], \ \ - \z^K(t)(0) =  \beta, \quad \forall t \in [0,T]
\end{equation}
and for all $S \in {\mathcal P}^+, \ T \in {\mathcal T}^+$, we have $\forall t \in
]0, T- \step]$
\begin{equation}\label{E2exist==}
 h(u^K(t+ \Delta t),  D_xT(u^K(t+ \Delta t)) \leq \z^K(t) D_xT(u^K(t+ \Delta t)) \ \ \
{\rm as \ measures}
\end{equation}
\begin{equation}\nonumber\label{E2exist1==}
\displaystyle h_S(u^K(t+ \Delta t),  D_xT(u^K(t+ \Delta t))
\leq \z^K(t) D_xJ_{T'S}(u^K(t+ \Delta t)) \ \ \ \hbox{as measures}.
\end{equation}
Note that (\ref{E2exist==}) is equivalent to
\begin{equation}\nonumber\label{E2exist}
 \frac{c}{2} \vert D_x^s((T(u^K(t+ \Delta t)))^2) \vert \leq \z^K(t) D_x^sT(u^K(t+ \Delta t)) \ \ \
{\rm as \ measures.}
\end{equation}

Since $\a(u^{n+1},(u^{n+1})') D_xT(u^{n+1}) \geq h(u^{n+1}, D_xT(u^{n+1}))$ as measures
in $]0,L[$,  using (\ref{BOUNFG}) we can write
$$h(u^{n+1}, D_xT(u^{n+1})) = \a(u^{n+1},(u^{n+1})')(T(u^{n+1}))' \mathcal{L}^1
+ \frac{c}{2} \vert D_x^s[(T(u^{n+1}))^2] \vert  $$
$$\geq \frac{c}{2} \vert ((T(u^{n+1}))^2)' \vert \mathcal{L}^1 - \frac{c^2}{\nu}  (T(u^{n+1}))^2 \mathcal{L}^1
+ \frac{c}{2} \vert D_x^s[(T(u^{n+1}))^2] \vert   $$
$$=  \frac{c}{2} \vert D_x[(T(u^{n+1}))^2] \vert - \frac{c^2}{\nu}  (T(u^{n+1}))^2 \mathcal{L}^1,$$
from  where we get the following inequality as measures
\begin{equation}\label{DesImport22}
\z^K(t) D_xT(u^K(t+ \Delta t)) \geq  \frac{c}{2} \vert D_x[(T(u^K(t+
\Delta t)))^2] \vert - \frac{c^2}{\nu} (T(u^K(t+ \Delta t)))^2.
\end{equation}

\begin{lemma} There exists $M:= M(\beta, c,\nu,L, \Vert u_0
\Vert_\infty)$ such that \label{small_dataK}
\begin{equation}\label{pestCL1}
\Vert u^K(t)\Vert_{\infty}\leq M \quad \quad \forall \, K \in \NN \
\hbox{and} \quad  \forall t \in [0,T].
\end{equation}
Consequently, $\Vert u(t)\Vert_{\infty}\leq M \quad \forall t \in [0,T].$
\end{lemma}

\noindent{\bf Proof.} Since
\begin{equation}\nonumber\label{Barsa2}
\left(I + \Delta t {\mathcal B }_{\beta} \right)^{-1}(u^n) =
u^{n+1}, \ \ \ \ \hbox{for} \ n=0,\ldots,K-1,
\end{equation}
 by Proposition \ref{Cojonudo99}, if $\mu:= \max \{ \frac{c \Vert
u_0 \Vert_{\infty}}{\beta}, 1 \}$, we have
$$0 \leq u^1 = \left(I + \Delta t {\mathcal B }_{\beta} \right)^{-1}(u_0)
\leq \mu
 u_{\beta}.$$
Then, repeating this process, we obtain
$$0 \leq  u^{n+1} = \left(I + \Delta t  {\mathcal B }_{\beta} \right)^{-1}(u^n)
\leq \left(I + \Delta t  {\mathcal B }_{\beta} \right)^{-1}(\mu
u_{\beta})
$$
$$= \mu \left(I + \Delta t  {\mathcal B
}_{\frac{\beta}{\mu}} \right)^{-1} (u_{\beta})\leq  \mu \left(I +
\Delta t  {\mathcal B }_{\beta} \right)^{-1}(u_{\beta}) = \mu
 u_{\beta},$$
and the proof concludes.

\hfill$\Box$

\noindent{\it Step 2.} By (\ref{pestCL1}), $\Vert \z^K(t)
\Vert_{\infty} \leq C$ \ for all $K \in \N$ and \ a.e. $t \in [0,
T]$. Then we may  assume that $\z^K \rightharpoonup \z \in L^{\infty}(Q_T)$ \  weakly$^*$.
Moreover, since $\z^K(t) = cu^K(t + \Delta t) \b (u^K(t + \Delta
t), \partial_x u^K(t + \Delta t)) \ \ \forall t \in ]0,T- \step]$,
with  $ \Vert \b (u^K(t + \Delta t), \partial_x u^K(t + \Delta t))
\Vert_\infty \leq 1$ and $u^K$ converges uniformly to $u$ in
$C([0,T],L^1(]0,L[))$, we may also assume that $\b (u^K(t + \step), \partial_x u^K(t + \Delta t))
\rightharpoonup \z_{\b}(t) \in L^{\infty}(Q_T) \
\hbox{weakly}^*$
and
\begin{equation}\label{EE2S2}
\z(t) = cu(t) \z_{\b}(t) \ \ \ \ \ {\rm for \ almost \ all} \ \ t
\in [0, T].
\end{equation}

 Given $w \in BV(]0, L[)$, from
(\ref{E1exist}) and (\ref{pestCL1}),
it follows that for each $t \in ]0,T]$
$$
\left\vert \int_0^L \xi^K(t,x) w(x) \, dx \right\vert = \left\vert
-\int_0^L \z^K(t)D w + \z^K(L)w(L) + \beta w (0_{+}) \right\vert
$$
$$
 \leq C \Vert w
\Vert_{BV(]0, L[)} + |\z^K(L)w(L)| \le (C+c \mu \|u_\beta\|_\infty) \Vert w
\Vert_{BV(]0, L[)},
$$
where the continuous injection of $BV(]0,
L[)$ into $L^\infty (]0, L[)$ was used. Thus, $\Vert\xi^K(t) \Vert_{BV(]0, L[)^*} \leq C, \
\forall \ K \in \N \  {\rm and} \  t \in ]0, T]$. Consequently,
$\{ \xi^K \}$ is a bounded sequence in \ $L^{\infty}(0, T; BV(]0,
L[)^*)$. Now, since  $L^{\infty}(0, T; BV(]0, L[)^*)$ is
a vector subspace of the dual space $\big(L^1(0, T; BV(]0,
L[))\big)^*$, we can find a subnet
 $\xi^\alpha$ of $\xi^K$ such that
\begin{equation}\nonumber
\label{E1S2} \xi^{\alpha} \rightharpoonup \xi \in \big(L^1(0, T;
BV(]0, L[))\big)^* \quad \hbox{weakly$^*$}.
\end{equation}

Working as in the proof of Theorem 5.5 of \cite{ARMA}, Step 2, we can prove that (\ref{IBPFTTT}) holds and $u_t = D_x \z \ \hbox{ in} \ \mathcal{D}^\prime (Q_T)$.

\noindent{\it Step 3.} Next, we prove that  $u_t = D_x\, \z$ in the sense given in (ii)
of Definition \ref{DefESevP}. To do this, let us first observe that we can prove, as in
the proof of Theorem 5.5 of \cite{ARMA}, Step 4, that the distribution $\z Dw$ in $Q_T$
defined by (\ref{DefMRAd}) is a Radon measure in $Q_T$ for all $w \in L^1(0, T;
BV(]0,L[))$, and also that
$$
  \langle \z D_x w , \varphi \rangle = \lim_\alpha \int_0^T \int_0^L \z^\alpha (t,x) D_x w(t,x) \varphi_x(t,x) \, dxdt.
$$
From where it follows, combining with (\ref{condel2existt}) and integrating by parts,
$$\int_{Q_T} \z D_xw  =  \lim_{\alpha} \int_0^T\int_0^L \z^\alpha(t) D_x w(t) \, dt = -\lim_{\alpha}\int_0^T \int_0^L w(t,x) D_x \z^\alpha (t,x) \,
dxdt
 $$ $$+ \lim_{\alpha} \left[\int_0^T \z^\alpha (t,L) w(t,L_{-}) \, dt - \int_0^T \z^\alpha(0) w(t,0_+) \, dt \right]$$
$$= \lim_{\alpha} \left[- \langle \xi^\alpha, w \rangle
  - c  \int_0^T u^{\alpha}(t+ \Delta t)(L_{-}) w(t,L_{-}) \, dt
  + \beta \int_0^T  w(t,0_+) \, dt \right] $$ $$= - \langle u_t, w \rangle
  - c  \int_0^T u(t)(L_{-}) w(t,L_{-}) \, dt+ \beta  \int_0^T w(t,0_+) \, dt,$$
and (\ref{IBPF}) holds.

\noindent{\it Step 4.} Let $T=T_{a,b}$ be any cut-off function, let $j$ be the primitive
of $T$. Let $0 \leq \phi \in {\mathcal D}(]0,T[)$. Multiplying (\ref{CL1}) by
$T(u^{n+1})\phi(t)$, $t\in (t_n,t_{n+1}]$ integrating in $(t_n,t_{n+1}]\times ]0,L[$ and
adding from $n=0$ to $n=K-1$, we have
\begin{equation}\label{Ec1CabreoA}\begin{array}{l}
    \displaystyle\sum_{n=0}^{K-1} \int_{t_n}^{t_{n+1}} \phi(t) \int_0^L \frac{u^{n+1}-u^n}{\step}  T(u^{n+1}) \, dxdt
 + \displaystyle \int_0^T \phi(t)\int_0^L \z^K(t) D_x( T(u^K(t+ \step))) \, dt \\ \\ = \displaystyle\sum_{n=0}^{K-1} \int_{t_n}^{t_{n+1}} \left(\beta \phi(t) T(u^{n+1}(0_+)) - c \phi(t) u^{n+1}(L_-) T(u^{n+1}(L_-))\right) \, dt.
  \end{array}
\end{equation}
Since $\phi$ has compact support in time in $(0,T)$, for $K$ large enough,  performing
like in (\ref{truco_estandar}), we  have
$$
- \sum_{n=0}^{K-1}  \int_{t_n}^{t_{n+1}} \phi(t) \int_0^L \frac{u^{n+1}-u^n}{\Delta t}
T(u^{n+1}) \, dxdt
 \leq  \int_0^T \int_0^L
j(u^K(t)) \frac{\phi(t)-\phi(t-\Delta t)}{\Delta t} \, dxdt.
$$
 Hence, from (\ref{Ec1CabreoA}) it follows that
\begin{equation}
\label{phi_constant}\begin{array}{ll}
  \displaystyle\int_0^T\int_0^L \z^K(t)  \phi(t) D_xT(u^K(t+ \step)) \, dt
\\ \\
 \le \displaystyle \int_0^T\int_0^L j(u^K(t)) \frac{\phi(t) - \phi(t-\step)}{\step}
\, dxdt  + \int_0^T  \beta \phi(t) T(u^K(t + \step,0_+)) \, dt.
\end{array}
\end{equation}
Given $\epsilon > 0$, if we take into (\ref{phi_constant}) any test $0 \le \phi \in
\mathcal{D}(]0,T[)$ such that $\phi(t)=1$ for $t \in ]\epsilon,T- \epsilon[$, having in
mind (\ref{pestCL1}), we get
\begin{equation}\nonumber\label{Nomelocc}\begin{array}{l}
 \displaystyle\int_\epsilon^{T-\epsilon}\int_0^L \z^K(t)  D_xT(u^K(t+ \step)) \, dt
\\ \\ \le \displaystyle \int_0^T\int_0^L j(u^K(t)) \frac{\phi(t) - \phi(t-\step)}{\step} \,
dxdt   +\displaystyle \int_0^T \beta T(u^K(t+ \Delta t,0_+)) \, dt \le  C.
\end{array}
\end{equation}
This implies that $\{ \z^K(t) D_x(T(u^K(t+\step))) \}$ is a bounded sequence in $\\
L_{loc,w}^1(0,T,\mathcal{M}(]0,L[))$, where $\mathcal{M}(]0,L[)$ denote the space of
bounded Radon measures in $]0,L[$.

On the other hand, by (\ref{DesImport22})

$$
   \int_{\epsilon}^{T-\epsilon} \int_0^L \z^K(t) D_x T(u^K (t+ \step)) \, dt
$$
$$
\ge \frac{c}{2} \int_{\epsilon}^{T-\epsilon} \int_0^L \vert D_x[(T(u^K(t+ \Delta t)))^2]
\vert \, dt - \int_{\epsilon}^{T-\epsilon} \int_0^L\frac{c^2}{\nu} (T(u^K(t+ \Delta
t)))^2\, dt.
$$
Hence
$$
    \int_{\epsilon}^{T-\epsilon} \int_0^L \vert
D_x[(T(u^K(t+ \Delta t)))^2] \vert \, dt \le
    \frac{2 C}{c} + \frac{2cLT b^2}{\nu} =  C,
$$
where by the coarea formula it follows that
\begin{equation}\label{AcotRRR}
    \int_{\epsilon}^{T-\epsilon} \int_0^L \vert
D_x T(u^K(t+ \Delta t)) \vert \, dt \le  C.
\end{equation}
Moreover, by Lemma 5 of \cite{ACM4:01}, the map $t \mapsto \Vert T(u^K(t))
\Vert_{BV(]0,L[)}$ is measurable, then by Fatou's Lemma and (\ref{AcotRRR}), it
follows that
\begin{equation}\label{AcotSDRE}\begin{array}{l}
    \displaystyle\int_{\epsilon}^{T- \epsilon} \liminf_{K \to \infty} \int_0^L |D_x T(u^K(t+ \step))| \, dt
\le  \displaystyle \liminf_{K \to \infty} \int_{\epsilon}^{T- \epsilon} \int_0^L
|D_xT(u^K(t+ \step))| \, dt \le C.
\end{array}
\end{equation}
Now, since the total variation is lower semi-continuous in $L^1(]0,L[)$, we have
$$\int_0^L \vert D_xT(u(t)) \vert \leq \liminf_{K \to \infty} \int_0^L \vert
D_xT(u^K(t)) \vert,$$ thus we deduce that $T(u(t)) \in BV(]0,L[)$ for almost all $t \in
(0,T)$ and consequently $u(t) \in TBV^+(]0,L[)$. Then, by (\ref{AcotSDRE}), applying
again Lemma 5 of \cite{ACM4:01}, we obtain that
\begin{equation}\label{loc1}
T(u(\cdot)) \in L^1_{loc, w}(0, T, BV(]0,L[)).
\end{equation}

\noindent{\it Step 5 (identification of the field).} Let us now prove that
\begin{equation}
\label{ident_field}
\z(t)= \a (u(t),\partial_x u(t)) \quad \mbox{a.e.} \, t \in ]0,T[.
\end{equation}

Let $0\le \phi \in \mathcal{D}(Q_T)$ and $g \in C^2([0,L])$. Assume that $\phi= \eta(t)
\rho(x)$ with $\eta \in \mathcal{D}(]0,T[)$ and $\rho \in \mathcal{D}(]0,L[)$. Let
$0<a<b$, and $T=T_{a,b}$.  Let $j$ denote the primitive of $T$. Recall that
$$
    J_{\a}(x,r) = \int_0^r \a(s,g^\prime(x)) \, ds \ \ \mbox{and} \quad
   J_{\a^\prime}(x,r)= \int_0^r \partial_x[\a(s,g^\prime(x))] \, ds
$$
For simplicity, we write
$$
D_2 J_{\a}(x,T(u^K(t+\step))):= D_x \left[J_{\a}(x,T(u^K(t + \step )))\right] -
J_{\a^\prime}(x,T(u^K(t + \step ))).
$$
Working as in the proof of step 6 of Theorem 3 in \cite{ARMA} we
find out that
\begin{equation}
\label{menosuno} [D_2 J_{\a}(x,T(u^K(t+\step)))]^{ac} = \a(u^K(t+\step), g')
\partial_x [T(u^K(t+\step))].
\end{equation}
Using (\ref{menosuno}), (\ref{monotoneineq}) and (\ref{E2exist}),
we obtain
$$
 \int_0^T \int_0^L \phi \ \z^K(t) D_x \left( T(u^K(t+\step)) - g \right) \,
dt
$$
$$
  - \int_0^T \int_0^L \phi \left[D_2 J_{\a}(x,T(u^K(t+\step))) -
\a(u^K(t+\step), g')g' \right] \, dt
$$
$$
   = \int_0^T \int_0^L \phi \left[\z^K(t) D_xT(u^K(t+ \step)) - \z^K(t)
g' + \a(u^K(t+\step),g') g' \right] \, dt
$$
$$
  -\int_0^T \int_0^L \phi \left\{[D_2 J_{\a}(x,T(u^K(t+\step)))]^{ac}+ [D_2
J_{\a}(x,T(u^K(t+\step)))]^s \right\} \, dt
$$
$$
 = \int_0^T \int_0^L \phi \left(\a(u^K(t+ \step),g') - \z^K(t)
\right) \left(g' - \partial_x T(u^K(t+ \step)) \right) \, dt
$$
$$
   + \int_0^T \int_0^L \phi \left[\z^K(t) D_x^sT(u^K(t+\step)) - [D_2 J_{\a}
(x,T(u^K(t+\step)))]^s \right] \, dt
$$
$$   \ge \int_0^T \int_0^L \phi \left[\z^K(t) D_x^sT(u^K(t +\step)))- [D_2
J_{\a} (x,T(u^K(t+\step)))]^s \right] \, dt
$$
$$   \ge \int_0^T \int_0^L \phi
\left[\frac{c}{2}\left|D_x^s(T(u^K(t+\step))^2)\right| - [D_2 J_{\a} (x,T(u^K(t+\step)))]^s
\right] \, dt.
$$
%
Again working as in the proof of step 6 of Theorem 3 in
\cite{ARMA}, we get
$$
    \int_0^T \int_0^L \phi \left[\frac{c}{2}\left|D_x^s(T(u^K(t+\step))^2)\right| -[D_2 (J_{\a}(x,T(u^K(t+ \step)))]^s \right] \, dt \ge 0.
$$
Therefore, we obtain
\begin{equation}
\label{eq_71}
\begin{array}{l}
\displaystyle \int_0^T \int_0^L \phi \ \z^K(t) D_x(T(u^K(t +
\step))-g) \, dt \\ \\ - \displaystyle \int_0^T \int_0^L \phi
\left[D_2 J_{\a}(x,T(u^K(t+ \step))) - \a(u^K(t+\step), g') g'
\right] \, dt \ge 0.
\end{array}
\end{equation}
Now we shall bound from above the first term. By (\ref{E1exist}) and for $\step$ small
enough,  performing like in (\ref{truco_estandar}), we get
\begin{equation}
\nonumber
\begin{array}{l}
\displaystyle
    \int_0^T\int_0^L \phi(t,x) T(u^K(t+\step)) D_x \z^K(t) \, dt
   = \int_0^T\int_0^L \phi(t,x) T(u^K(t+\step) \xi^K(t) \, dxdt
\\ \\
\displaystyle
   \ge \int_0^T \int_0^L \frac{\phi(t- \step,x) - \phi(t,x)}{\step} j(u^K(t)) \, dtdx
    \end{array}
\end{equation}
Then,  integrating by parts, we have
$$
    \int_0^T \int_0^L \phi(t) \ \z^K(t) D_x(T(u^K(t+ \step))-g) \, dt \le - \int_0^T \int_0^L \frac{\phi(t- \step) - \phi(t)}{\step}
j(u^K(t)) \, dtdx
$$
$$
   + \int_0^T \int_0^L \phi(t) g \xi^K(t) \, dtdx - \int_0^T \int_0^L \partial_x \phi(t) \
   \z^K(t) [T(u^K(t+\step))-g]\, dtdx.
$$
Thanks to this inequality we arrive from (\ref{eq_71}) to
\begin{equation}
\label{eq_72}
\begin{array}{l}
\displaystyle - \int_0^T \int_0^L \frac{\phi(t- \step) - \phi(t)}{\step} j(u^K(t)) \, dtdx + \int_0^T \int_0^L \phi(t) g \xi^K(t) \, dtdx
\\ \\
\displaystyle - \int_0^T \int_0^L \partial_x \phi(t) \ \z^K(t) [T(u^K(t+\step))-g]\, dtdx
\\ \\
- \displaystyle \int_0^T \int_0^L \phi(t) \left[D_2 J_{\a}(x,T(u^K(t+ \step))) -
\a(u^K(t+\step), g') g' \right] \, dt \ge 0.
\end{array}
\end{equation}

Letting $K \to \infty$ in (\ref{eq_72}) and having in mind that
$$D_2 J_{\a}(x,T(u^K(t + \step))) \rightharpoonup D_2
J_{\a}(x,T(u(t))) \quad \hbox{ weakly as measures}$$ we obtain
\begin{equation}
\label{eq_72.5}
\begin{array}{l}
\displaystyle \int_0^T \int_0^L \partial_t \phi(t) j(u(t)) \, dt +
\langle u_t,\phi g \rangle \, - \int_0^T \int_0^L [T(u(t))-g]\z(t) \partial_x \phi(t) \,
dxdt
\\ \\
+\displaystyle \int_0^T \int_0^L \phi(t) \left[-D_2 J_{\a}(x,T(u(t)))+ \a(u(t),g')g'\right] \, dt \ge 0.
\end{array}
\end{equation}
By (\ref{IBPF}),
$$
     \langle u_t, \phi \ g  \rangle
     = - \int_0^T \int_0^L \z(t) g \partial_x \phi(t) \, dtdx -
\int_0^T \int_0^L  \z(t) g' \phi(t) \, dtdx
$$
and we can rearrange (\ref{eq_72.5}) in the following way
\begin{equation}
\label{eq_73}
\begin{array}{l}
   \displaystyle \int_0^T \int_0^L \partial_t \phi(t) j(u(t)) \, dtdx - \int_0^T \int_0^L  \z(t) g'
\phi(t) \, dtdx \displaystyle   - \int_0^T \int_0^L T(u(t)) \z(t)
\partial_x \phi(t) \, dxdt
\\ \\
\displaystyle   + \int_0^T \int_0^L \phi(t) \left[-D_2 J_{\a}(x,T(u(t)))+ \a(u(t),g')g'\right] \, dt \ge 0.
\end{array}
\end{equation}

Now, for $\tau$ small enough and using again the trick in (\ref{truco_estandar}), we have
$$
   \int_0^T \int_0^L \partial_t \phi(t,x) j(u(t,x)) \, dxdt =\lim_{\tau \to 0} \int_0^T \int_0^L
\frac{\eta(t - \tau) - \eta(t)}{-\tau} j(u(t,x)) \rho(x) \, dxdt
$$
$$
\displaystyle \leq \lim_{\tau \to 0} \int_0^T \int_0^L u(t,x) \rho(x) \frac{d}{dt} (\eta
T(u))^\tau(t,x) \, dxdt,
$$
where we used again the notion of Dunford integral (see Remark \ref{timereg}). Using
(\ref{IBPFTTT}), we have
$$
   \int_0^T \int_0^L u(t) \rho \frac{d}{dt} (\eta T(u))^\tau (t) \, dxdt
   = - \langle u_t,\rho (\eta T(u))^\tau(\cdot) \rangle
   $$
$$
   = - \lim_{\alpha}  \langle \xi^\alpha,\rho(\eta T(u))^\tau(\cdot)\rangle
    = - \lim_{\alpha} \int_0^T
  \langle D_x\z^\alpha(t), \rho \frac{1}{\tau} \int_{t-\tau}^t \eta(s) T(u(s)) \, ds\rangle \, dt
$$
$$
  = \lim_{\alpha} \int_0^T \int_0^L \z^\alpha(t) D_x \left(\rho
\frac{1}{\tau} \int_{t - \tau}^t \eta(s) T(u(s)) \, ds \right) \,
dt
   =  \lim_{\alpha} \int_0^T \int_0^L \partial_x \rho \ \z^\alpha(t)
   \int_{t-\tau}^t \frac{1}{\tau} \eta(s) T(u(s))
   \, dsdxdt
  $$
  $$
+ \lim_{\alpha} \int_0^T \int_0^L\rho \
\z^\alpha(t)D_x[(\eta T(u))^\tau (t)] \, dt
   = \int_0^T \frac{1}{\tau} \int_{t-\tau}^t \eta(s) \int_0^L T(u(s)) \z(t) \partial_x \rho
   \, dxdsdt
$$
$$
 +\lim_{\alpha} \int_0^T \int_0^L \rho \ \z^\alpha(t)
\partial_x[(\eta T(u))^\tau (t)] \, dxdt
 + \lim_{\alpha} \int_0^T \int_0^L \rho \ \z^\alpha(t)
D^s_x[(\eta T(u))^\tau (t)] \,dt
$$
$$ \le \int_0^T \frac{1}{\tau} \int_{t-\tau}^t \eta(s) \int_0^L T(u(s)) \z(t)
\partial_x \rho \, dxdsdt
 + \int_0^T \frac{1}{\tau} \int_{t -\tau}^t \eta(s) \int_0^L \rho \z(t)
\partial_x(T(u(s))) \, dxdsdt
$$
$$
 + \int_0^T \frac{1}{\tau} \int_{t-\tau}^t \eta(s) \int_0^L cM \rho \vert
D^s_x[T(u(s)] \vert \, ds dt.
$$
Taking limits when $\tau \to 0$, having in mind (\ref{pestCL1}), we obtain
$$
   \int_0^T \int_0^L \partial_t \phi(t) j(u(t)) \, dxdt \le \int_0^T \eta(t) \int_0^L T(u(t)) \z(t)
\partial_x \rho \, dxdt
$$
$$
   + \int_0^T \eta(t) \int_0^L \rho \ \z(t) \partial_x T(u(t)) \, dxdt
+  c M \int_0^T \eta(t) \int_0^L \rho  \vert D^s_x[T(u(t)] \vert \, dt.
$$
From (\ref{eq_73}), all gathered together reads
$$
   0 \le - \int_0^T \int_0^L \phi(t) \z(t) g' \, dxdt + \int_0^T \eta(t)
   \int_0^L \rho \ \z(t) \partial_x (T(u(t))) \, dxdt
    + c M \int_0^T   \eta (t) \int_0^L \rho   \vert D^s_x T(u(t)) \vert  \, dt
    $$ $$ + \int_0^T \int_0^L \phi \left[-D_2 J_{\a}(x,T(u(t)))+ \a(u(t),g')g'\right]
    \, dt.
$$
Using that $  D_2 J_{\a} (x,T(u(t)) = \a(u(t),g') \partial_x (T(u(t))) + [D_2
J_{\a}(x,T(u(t)))]^s$, this is written as
\begin{equation}
\nonumber
\begin{array}{l}
\displaystyle
  0 \le c M \int_0^T   \eta (t) \int_0^L \rho   \vert D^s_x T(u(t)) \vert \, dt  -
  \int_0^T \int_0^L \phi [D_2J_{\a}(x,T(u(t))]^s \, dt
\\
\displaystyle
   + \int_0^T \int_0^L [g' - \partial_x (T(u(t)))][\a(u(t),g') -
   \z(t)] \phi \, dxdt.
\end{array}
\end{equation}
As measures,
$$
  cM \vert D^s_x T(u(t)) \vert - [D_2J_{\a}(x,T(u(t))]^s + [g' -
  \partial_x (T(u(t)))][\a(u(t),g') - \z(t)] \mathcal{L}^2 \ge 0.
$$
Taking the absolutely continuous part and particularizing to points
$x \in [a<u(t)<b]$, this reduces to
 $$
  [g' - \partial_x u(t)][\a(u(t),g') - \z(t)] \ge 0,
  $$
   an inequality which holds for all $(t,x) \in S \cap [a < u < b]$, where $S \subseteq ]0,T[\times ]0,L[$ is such that $\mathcal{L}^2(]0,T[ \times ]0,L[ \setminus S) = 0$, and all $g \in C^2([0,L]).$
Being $(t,x) \in S \cap [a < u < b]$ fixed and $\xi \in \RR$ given, we can find a function $g$ as above such that $g' (x) = \xi$. Then
$$
    \left(\z(t,x) - \a(u(t),\xi)\right) (\partial_x u(t,x) - \xi) \ge 0, \, \, \forall \xi \in \RR \quad \mbox{and}\, \forall (t,x) \in S\cap [a < u <b].
$$
By an application of Minty-Browder's method in $\RR$, these inequalities imply that
$$
   \z(x) = \a(u(t,x), \partial_x u (t,x)) \quad \mbox{a.e. on}\, Q_T \cap [a < u <b].
$$
Since this holds for any $0<a<b$, we obtain (\ref{ident_field}) a.e. on the points of
$Q_T$ such that $u(t,x) \ne 0$. Now, by our assumptions on $\a$ and (\ref{EE2S2}) we
deduce that $\z(x) = \a(u(x),u^\prime(x)) = 0$ a.e. on $[u=0]$. We have proved
(\ref{ident_field}).

\noindent{\it Step 6. The entropy inequality.} Given $S\in {\mathcal P}^+, T \in \mathcal{T}^+$ and $\phi
\in \mathcal{D}(Q_T)$, working as in the proof of (\ref{phi_constant}) we can get
\begin{equation}
\label{phi_out_compuesta}
\begin{array}{l}
\displaystyle
   \int_0^T \int_0^L \phi \ \z^K(t) D_x \left(T(u^K(t+ \step)) S(u^K(t+ \step)) \right) \, dt
\\ \\
\displaystyle
   \le \int_0^T \int_0^L J_{TS}(u^K(t)) \frac{\phi(t) - \phi(t - \step)}{\step} \, dxdt
\\ \\ \displaystyle
   - \int_0^T \int_0^L \z^K(t) \partial_x \phi T(u^K(t+ \step)) S(u^K(t+ \step)) \, dxdt
 \end{array}
\end{equation}
and the fact that $  \left\{\z^K(t) D_x \left(T(u^K(t+ \step)) S(u^K(t+ \step)) \right) \right\} $
 is a bounded sequence in $L_{loc}^1(0,T;\mathcal{M}(]0,L[))$. From here, as in the proof
of Theorem 4.5 of \cite{ACMRelat}, we can get that the sequences
$\{\z^K(t) D_xJ_{T'S}(u^K(t+ \step)) \}$ and $\{\z^K(t)
D_xJ_{S'T}(u^K(t+ \step)) \}$ are bounded in
$L_{loc}^1(0,T;\mathcal{M}(]0,L[))$. This allows us to define, up
to subsequence, the objects $\mu_T^S, \mu_S^T\in \mathcal{M}(Q_T)$
by means of
$$
  \langle \phi, \mu_S^T\rangle = \lim_{K} \int_0^T \int_0^L \phi \ \z^K(t) D_xJ_{T'S}(u^K(t+ \step)) \, dt, \quad \forall \phi \in C_c(Q_T),
$$
$$
 \langle \phi, \mu_T^S\rangle  = \lim_{K} \int_0^T \int_0^L \phi \ \z^K(t) D_xJ_{S'T}(u^K(t+ \step)) \, dt , \quad \forall \phi \in C_c(Q_T).
$$
Then, passing to the limit in (\ref{phi_out_compuesta}), we obtain
\begin{equation}\label{EccasiFF}\begin{array}{l}
    \displaystyle\langle \phi, \mu_S^T \rangle +\langle\phi, \mu_T^S \rangle \le
    \int_0^T \int_0^L J_{TS}(u(t)) \partial_t \phi(t) \, dxdt
\\ \\
   - \displaystyle \int_0^T \int_0^L \z(t) \partial_x \phi T(u(t))S(u(t)) \, dxdt, \quad \forall \phi \in \mathcal{D}(Q_T),
\end{array}
\end{equation}
Working as in proof  of Lemma 4.11 in \cite{ACMRelat}, we can get the following result.
\begin{lemma}
\label{lema4.11}
For $S\in {\mathcal P}^+, T \in \mathcal{T}^+$, we have that $\mu_S^T \ge h_S(u,DT(u))$.
\end{lemma}

By the above lemma and (\ref{EccasiFF}) we obtain the entropy inequality
$$\displaystyle \int_0^T\int_0^L \phi
h_S(u,DT(u)) \, dt + \int_0^T\int_0^L \phi h_T(u,DS(u)) \, dt
 \leq$$ $$ \displaystyle\int_0^T\int_0^L J_{TS}(u) \phi^{\prime} \, dxdt-
\int_0^T\int_0^L \a(u, \partial_x u) \cdot \partial_x \phi \ T(u) S(u) \, dxdt$$
 for truncatures $S\in {\mathcal P}^+ , \, T \in {\mathcal T}^+$ and any  smooth function $\phi$ of
 compact support.

\hfill$\Box$

\medskip \noindent {\it Acknowledgements.}  The first and
third  authors have been partially supported by the Spanish MCI
and FEDER, project MTM2008-03176. The second and fourth authors
have been partially supported by the Spanish MCI and FEDER,
project MTM2008-05271 and Junta de Andaluc\'{\i}a FQM--4267.

\end{document}